\documentclass[conference]{IEEEtran}
\pdfoutput=1
\usepackage{color}
\usepackage{graphicx}
\usepackage{enumerate}
\usepackage{amsmath,amssymb}
\usepackage{caption}
\usepackage{subcaption}

\newtheorem{theorem}{Theorem}
\newtheorem{definition}{Definition}

\DeclareMathOperator{\Proj}{Proj}

\DeclareMathOperator{\Int}{int}
\DeclareMathOperator{\co}{co}

\newcommand{\RE}{\mathbb{R}}
\newcommand{\argmin}{\operatornamewithlimits{argmin}}

\title{Implementing robust neuromodulation in neuromorphic circuits}
\author{
 \IEEEauthorblockN{Fernando Casta\~nos$^\dag$}\IEEEauthorblockA{Centro de Investigaci\'on y de Estudios Avanzados \\ 
  Instituto Polit\'ecnico Nacional, M\'exico \\ Email: castanos@ieee.org}\\
$^\dag$ {\small The two authors contributed equally to this work}
\and
 \IEEEauthorblockN{Alessio Franci$^\dag$}\IEEEauthorblockA{Universidad National Autonoma de M\'exico 
  \\ Facultad de Ciencias \\ Email: afranci@ciencias.unam.mx}

}

\date{\today}

\begin{document}

\maketitle

\begin{abstract}
We introduce a methodology to implement the physiological transition {between distinct neuronal spiking modes} in electronic circuits composed of resistors, capacitors and transistors. The result is a simple neuromorphic device organized by the same geometry {and exhibiting the same input--output properties as} high-dimensional electrophysiological neuron models. {Preliminary} experimental results highlight the robustness of the approach in real-world applications.
\end{abstract}

\section{Introduction}

Nature offers spectacular examples of energy-efficient, lightweight control
architectures. Flight control in an animal like a honeybee outperforms
the latest robotic architectures in terms of energy consumption, adaptability,
robustness, and dimensions. Neuromorphic engineering aims at emulating
the way in which biological neuronal systems perceive and represent the
outside world, take decisions and develop computations, and command motor outputs \cite{Boahen2005,Liu2010}.

In implementing the dynamical behavior of biological neurons in electronic hardware
we face the compromise between fidelity of the reproduced behavior and complexity
of the designed circuit. Existing silicon neuron designs span a variety of solutions:
from detailed implementation of neuron biophysics \cite{Hynna2007} to implementation of simple, abstract
neuron models \cite{Wijekoon2008}. Both approaches have advantages and disadvantages, and
it is an active research area to determine which implementation to use depending on
the desired objective \cite{Indiveri2011}.

The possibility of robustly and rapidly switching between distinct activity modes is one
of the peculiarity of biological neurons, which allows them to adapt their input--output
response to internal and environmental conditions. Two fundamental neuronal
activity modes are tonic spiking and bursting. Tonic spiking describes the slow, regular
generation of spikes in the neuron membrane potential. Bursting describes the alternation
between moments of low membrane potential and moments of high oscillatory activity, in
which spikes are generated at very high frequency. The transition between tonic spiking
and bursting plays a major role in {neuronal information processing by modulating neuron input--output behavior} \cite{Krahe2004,Sherman2001,lee2012neuromodulation}.

We showed~\cite{Franci2013,Franci2014} that all biological neurons share the same geometry
at the transition between {distinct spiking modes}. In particular, this transition can be
described in a simple, abstract model given by the normal form of an organizing singularity. Roughly speaking,
a singularity describes a highly degenerate and fragile situation that correspond to the transition 
between distinct regimes~\cite{Schaeffer1985}. There is a direct correspondence between
biophysical parameters and mathematical parameters in the abstract model, which leads
to a novel mathematical understanding of robustness and modulation of neuronal activity \cite{Drion2015}. 
We further showed that the same qualitative picture can be realized in simple {input--output} circuits \cite{Franci2014a}.

In this paper we follow the recipe provided in~\cite{Franci2014a} to design a neuromorphic
circuit with the property of exhibiting the same qualitative geometry, robustness, modulation capabilities, {and input--output behavior} as biophysical neuron models. As a first, biologically relevant illustration, we focus on the transition between tonic spiking and bursting.
The resulting circuit solely uses six transistors and passive elements {and its robust real-world implementation in low-cost components solely uses four additional transistors to overcome loading effects.

The key contributions with respect to existing neuromorphic design methods \cite{Indiveri2011,misra2010,almasi2016} are threefold. First, the equivalence (from a geometric, dynamical systems, and input--output perspective) of the designed circuit and high-dimensional biophysical neuron models close to the transition between distinct spiking modes is a provable consequence of the used approach. Second, the geometric nature of our methodology avoids laborious and non-constructive parameter fitting procedures, and, third, it also ensures robustness to components variability in real-world applications. }

In Section \ref{SEC: geoemetry of neuronal bursting} we rapidly review the results in~\cite{Franci2013,Franci2014,Franci2014a}.
Grounded in these works, we derive an implementation of our neuromorphic device in Section \ref{SEC: electronica} and
simulate it in \texttt{ngspice}~\cite{ngspice} (the code can be found in the Appendix). The actual implementation of this circuit and preliminary experimental tests are reported in Section~\ref{SEC: experimental}. {Future directions are discussed in Section~\ref{SEC: conclusions}.}

\section{The geometry of neuronal behaviors and its block realization} \label{SEC: geoemetry of neuronal bursting}

\begin{figure}
\centering
\includegraphics[width=0.90\columnwidth]{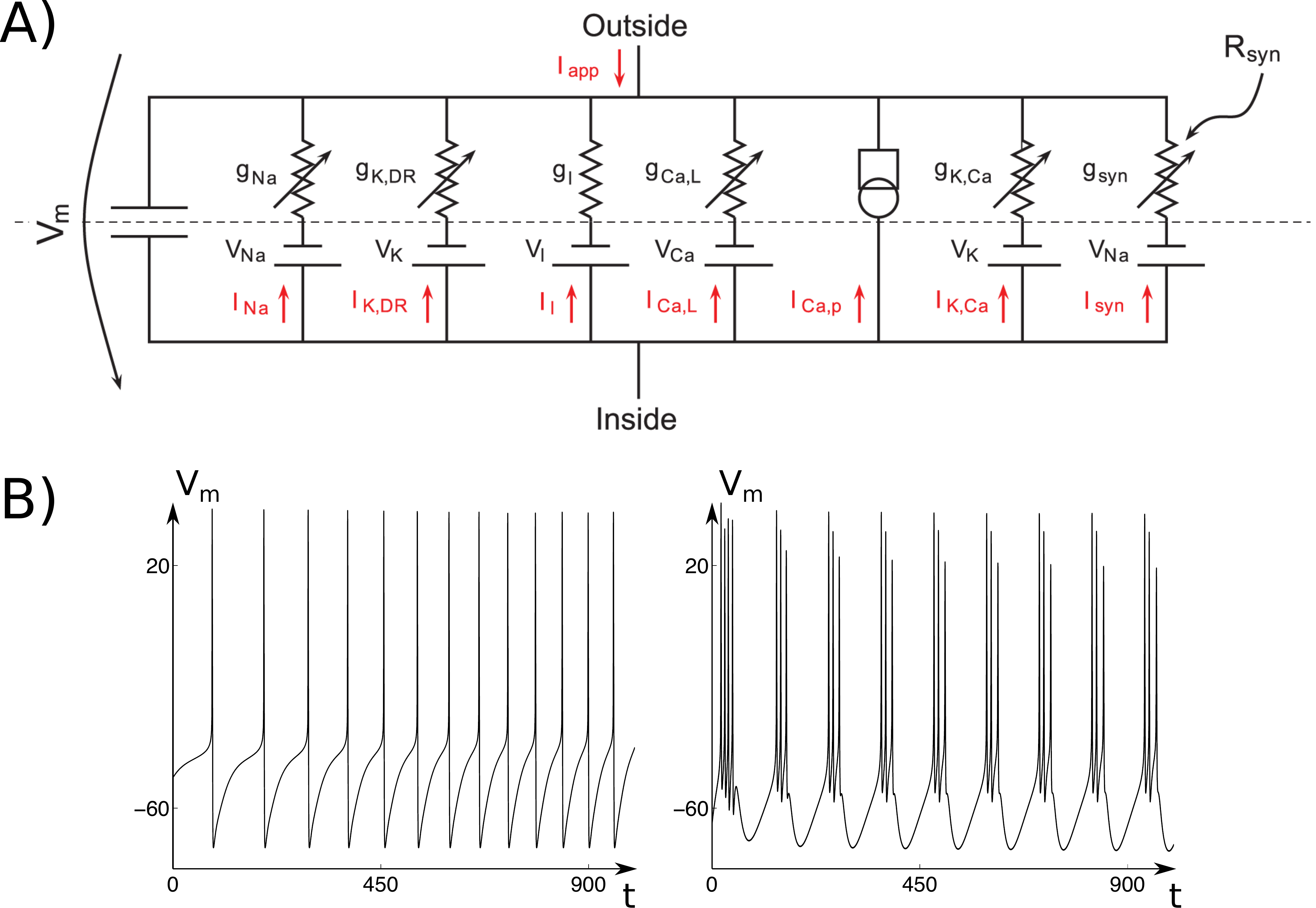}
\caption{A. The RC circuit associated to a conductance-based biophysical model of dopaminergic neuron (adapted from \cite{Drion2011}). B. Two activity modes in a conductance-based model.  Left: tonic spiking. Right: bursting. Membrane potential values are in millivolts. Time is in millisecond.}\label{FIG: neuron RC circuit}
\end{figure}

Electrophysiological models of neurons are constructed upon the seminal work of Hodgkin and Huxley \cite{Hodgkin1952}. They all share the physical interpretation of the nonlinear RC circuit depicted in Fig.~\ref{FIG: neuron RC circuit}-A. The capacitor models the neuron lipidic membrane and the other branches, containing a voltage source and a variable resistance, model the flow of a specific ion across the membrane.

Ion flow across the membrane is dynamically regulated by the membrane potential via opening and closing of {\it ion channels}, which makes the circuit in Fig.~\ref{FIG: neuron RC circuit}-A highly nonlinear. As such, it can exhibit a rich variety of dynamical behaviors. The present paper focuses on two fundamentals behaviors shared by almost all neuron types: the tonic spiking behavior of Fig.~\ref{FIG: neuron RC circuit}-B left  and the bursting behavior of Fig.~\ref{FIG: neuron RC circuit}-B right.

Reproducing tonic spiking and bursting, as well as the transition between these two modes, in an electrophysiological model requires an empirical tuning of the many biophysical parameters that usually ends up in an extensive brute-force computational parameter search \cite{Prinz2003}. A different approach relies on bifurcation theory \cite{Guckenheimer1983}.

Roughly speaking, bifurcation theory makes the ansatz that the vector field associated to an electrophysiological model undergoes some qualitative change at the transition between two distinct activity modes.

We showed in \cite{Franci2013,Franci2014} that the bifurcation associated to the transition between tonic spiking and bursting can be algebraically tracked by exploiting the multi-timescale nature of electrophysiological neuron models and by detecting a transcritical singularity in the critical manifold of the associated singularly perturbed dynamics. We refer the reader to \cite{Jones1995} for an introduction to geometric singular perturbation theory and to \cite{Schaeffer1985} for singularity theory concepts.

The power of this analysis is that we can visualize the geometry of the tonic spiking--bursting transition in a low-dimensional normal form of the organizing singularity:
\begin{subequations}\label{EQ:bursting normal form}
\begin{align}
 \dot x &= -x^3-(\lambda+y)^2+\beta x-(\alpha+u)-z \\
 \dot y &= \varepsilon_s(x- y) \\
 \dot z &= \varepsilon_u(x-z) \;, 
\end{align}
\end{subequations}
where $\lambda$ is called the bifurcation parameter, $\alpha,\beta$ are called unfolding parameters, and $0<\varepsilon_u\ll\varepsilon_s\ll 1$ model timescale separation between the three state variables $x$, $y$ and $z$. The system is driven by the external input $u$. The distinction between bifurcation and unfolding parameters is instrumental to the tools used in the construction of the normal form (\ref{EQ:bursting normal form}), that is, singularity theory applied to bifurcation problems \cite{Schaeffer1985}.

{The main theorem in \cite{Franci2014} provides constructive conditions on the parameters in model (\ref{EQ:bursting normal form}) to enforce the existence of a bursting attractor. This attractor exists for parameters close to the same transcritical singularity organizing the transition between tonic spiking and bursting in biophysical conductance-based models. As a corollary, the results in \cite{Franci2014} therefore provide a geometric way of exploring the transition between bursting and tonic spiking.
}

\begin{figure}
\centering
\includegraphics[width=0.99\columnwidth]{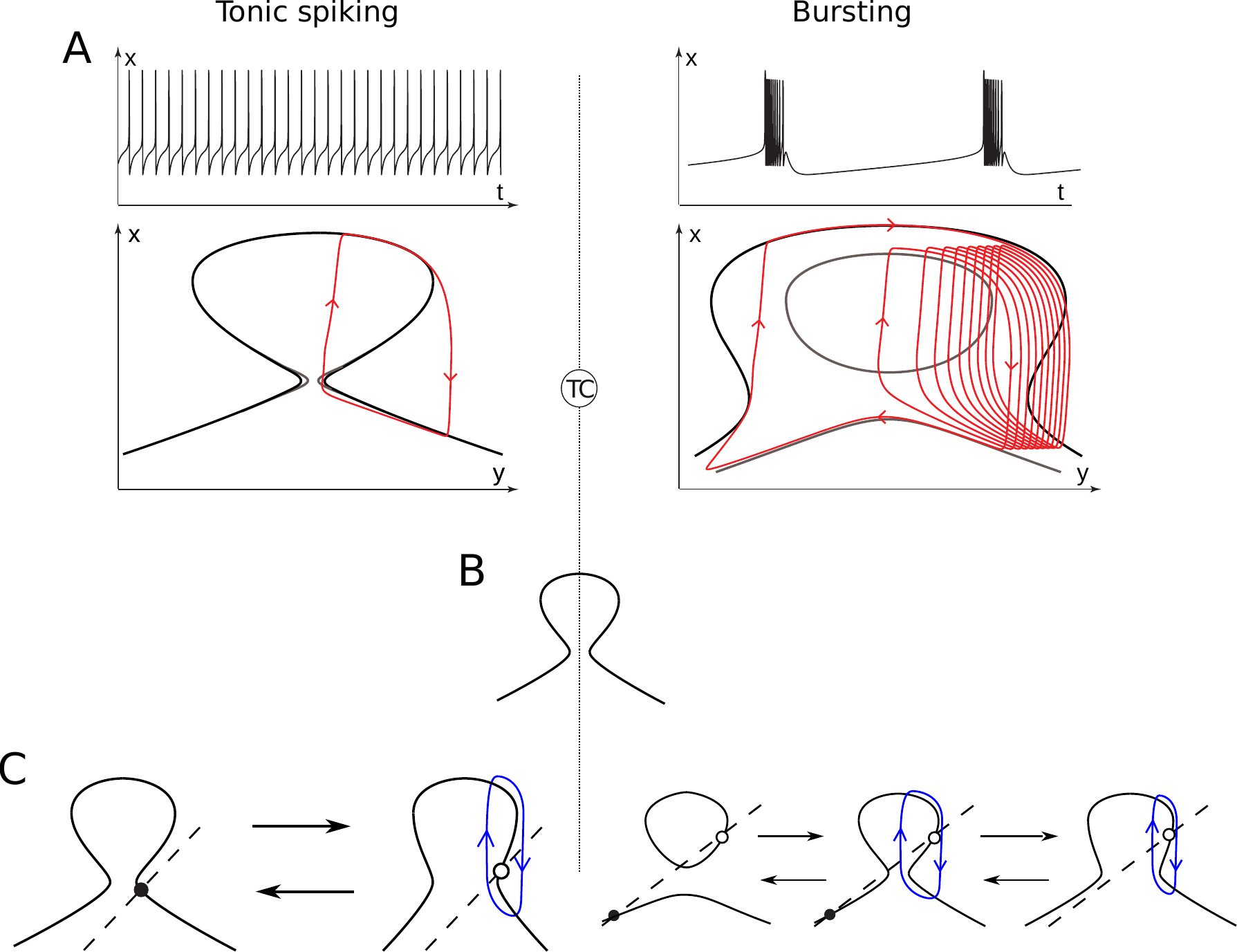}
\caption{A. Top: temporal traces of mode (\ref{EQ:bursting normal form}) in tonic spiking and bursting. Bottom: projection of the trajectory on the $x-y$ phase plane. The trajectory is in red. The critical manifold (\ref{EQ:bursting normal form critical manifold}) is drawn in black for a small value of $z$ (that is, at the beginning of spikes/bursts) and in gray for a large value of $z$ (that is, when the trajectory relaxes to rest). B. The mirrored hysteresis bifurcation diagram. C. The geometry of tonic spiking and bursting. The black dot denotes stable fixed point, the circle unstable fixed points. Stable limit cycles are drawn in blue.}\label{FIG:bursting phase portrait}
\end{figure}

Fig.~\ref{FIG:bursting phase portrait}-A shows the temporal traces and the projection onto the phase plane of the slow--fast subsystem (\ref{EQ:bursting normal form}a)-(\ref{EQ:bursting normal form}b) of tonic spiking and bursting behaviors in model (\ref{EQ:bursting normal form}). Due to timescale separation, trajectories spend most of the time near the critical manifold
\begin{equation}\label{EQ:bursting normal form critical manifold}
\mathcal Z:=\{(x,y,z)\in\mathbb R^3:\ -x^3-(\lambda+y)^2+\beta x-\alpha-z=0\} \;,
\end{equation}
that is, the $x$-nullcline composed of $x$ steady states as $y$ and $z$ vary.

The fundamental {geometric shape underlying the transition between tonic spiking and bursting is the mirrored hysteresis bifurcation diagram} introduced in \cite{Franci2014} and sketched in Fig.~\ref{FIG:bursting phase portrait}-B. In the tonic spiking mode, trajectories solely visit the right branch of the mirrored hysteresis {organizing the critical manifold}, whereas in bursting  mode trajectories alternate between the two branches.

The geometry of both behaviors is {summarized} in Fig.~\ref{FIG:bursting phase portrait}-C. In tonic mode, for a fixed value of the ultra-slow variable $z$, the slow--fast subsystem possesses a single attractor on the right branch of the mirrored hysteresis: either a stable fixed point, for a large value of $z$ (left plot), or a stable limit cycle, for a small value of $z$ (right plot). For large initial values of $z$ the model is therefore at quasi-steady state. In this case the ultra-slow dynamics (\ref{EQ:bursting normal form}c) forces $z$ to decrease until the steady state looses stability and a spike is emitted along the newborn limit cycle. This in turn leads to a sharp increase of $z$, which immediately lets the steady state recover stability, and the model is forced back to quasi-steady state {after a single spike}. 

In bursting mode, there exists a large range of values of the ultra-slow variable $z$ in which the slow--fast subsystem exhibits \emph{bistability} between a stable steady state on the left branch of the mirrored hysteresis and a limit cycle on the right branch of the mirrored hysteresis (center plot). Bursting arises from ultra-slow hysteretic evolution around this bistable region.

For a large initial value of $z$ the sole attractor is a stable steady state (left plot). The ultra-slow dynamics (\ref{EQ:bursting normal form}c) forces $z$ to decrease. However, once in the bistable region, the model remains at quasi-steady state. Only for sufficiently small $z$ the steady state looses stability and the trajectory converges toward the spiking limit cycle, which is now the sole attractor (left plot). On the limit cycle the ultra-slow dynamics (\ref{EQ:bursting normal form}c) forces $z$ to {increase}. Again, all through the bistable region the system remains  in the oscillatory mode and only for $z$ sufficiently large the trajectory converges back to the quasi-steady state.


Normal forms are useful not only because they unmask the geometry underlying a given dynamical behavior, but also because they possess the minimum number of parameters to reproduce a family of behaviors of interest. As such they can also be implemented in physical devices more easily and robustly than the original biophysical model, yet, preserving the same geometric and input--output properties.

\begin{figure}
\centering
\includegraphics[width=0.90\columnwidth]{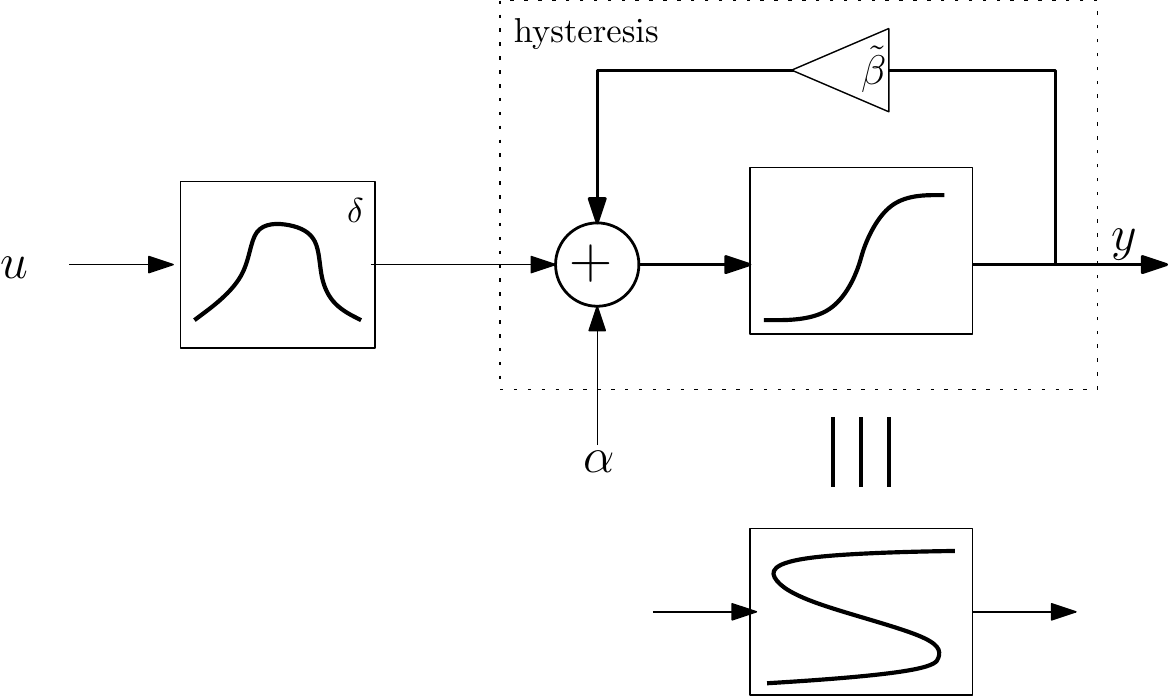}
\caption{Block realization of a mirrored hysteresis nonlinear characteristic (adapted from \cite{Franci2014a}).}\label{FIG: wcusp circuit}
\end{figure}

We showed in \cite{Franci2014a} that the mirrored hysteresis bifurcation diagram at the core of model (\ref{EQ:bursting normal form}) can be realized in the input--output diagram of Fig.~\ref{FIG: wcusp circuit}. Its basic ingredients are a non-monotone nonlinearity cascaded with a saturation nonlinearity and a positive feedback loop around the saturation nonlinearity. The positive feedback loop transforms the saturation into a hysteresis \cite[Proposition 1]{Franci2014a}, the non-monotone block creates a mirror of the resulting hysteretic characteristic.

\begin{figure}
\centering
\includegraphics[width=0.90\columnwidth]{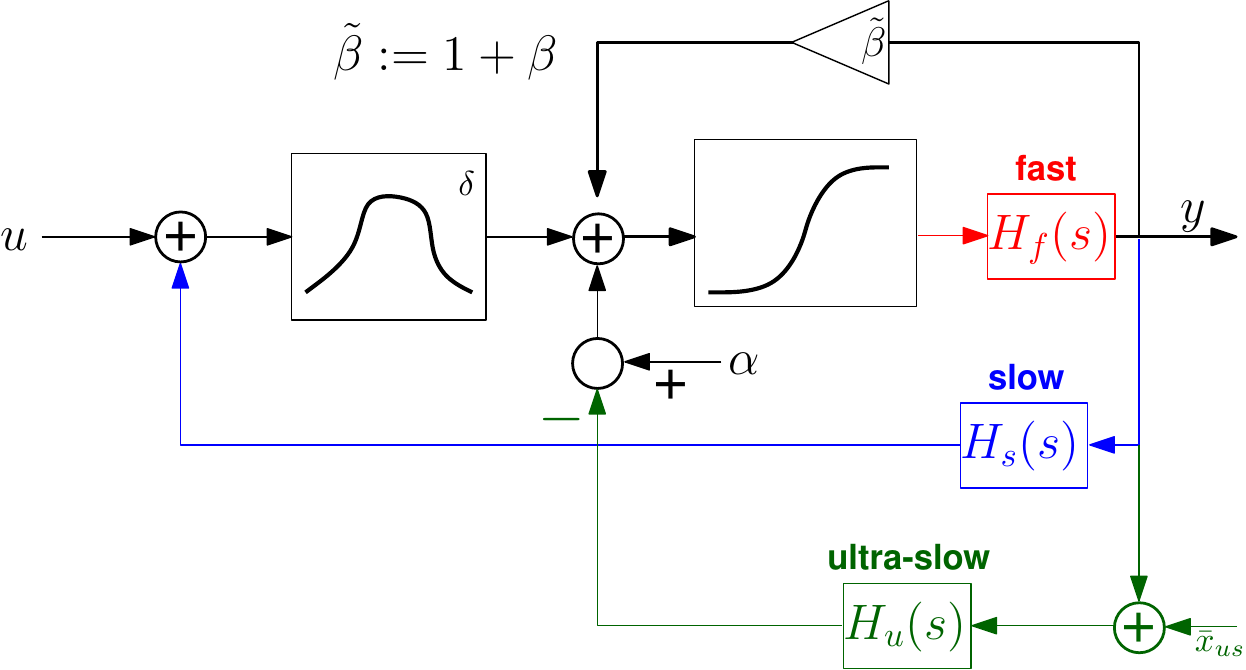}
\caption{Block realization of a neuron model (adapted from \cite{Franci2014a}). The linear filters $H_{f,s,u}(s)$ are first order filters with sharply separated cut-off frequency: large for $H_f$, intermediate for $H_s$, and small for $H_u$.}\label{FIG: burst circuit}
\end{figure}

Adding linear dynamical systems evolving on three sharply different timescales (Fig.~\ref{FIG: burst circuit}) transforms the static diagram in Fig.~\ref{FIG: wcusp circuit} into a three-timescale dynamical system with the same qualitative behavior of model (\ref{EQ:bursting normal form}). In particular, the diagram in Fig.~\ref{FIG: burst circuit} exhibits the same geometric transition between tonic spiking and bursting as model (\ref{EQ:bursting normal form}) \cite[Theorem 3]{Franci2014a}.

\section{Electronic implementation} \label{SEC: electronica}

\begin{figure}
\centering
 \includegraphics[scale=0.29]{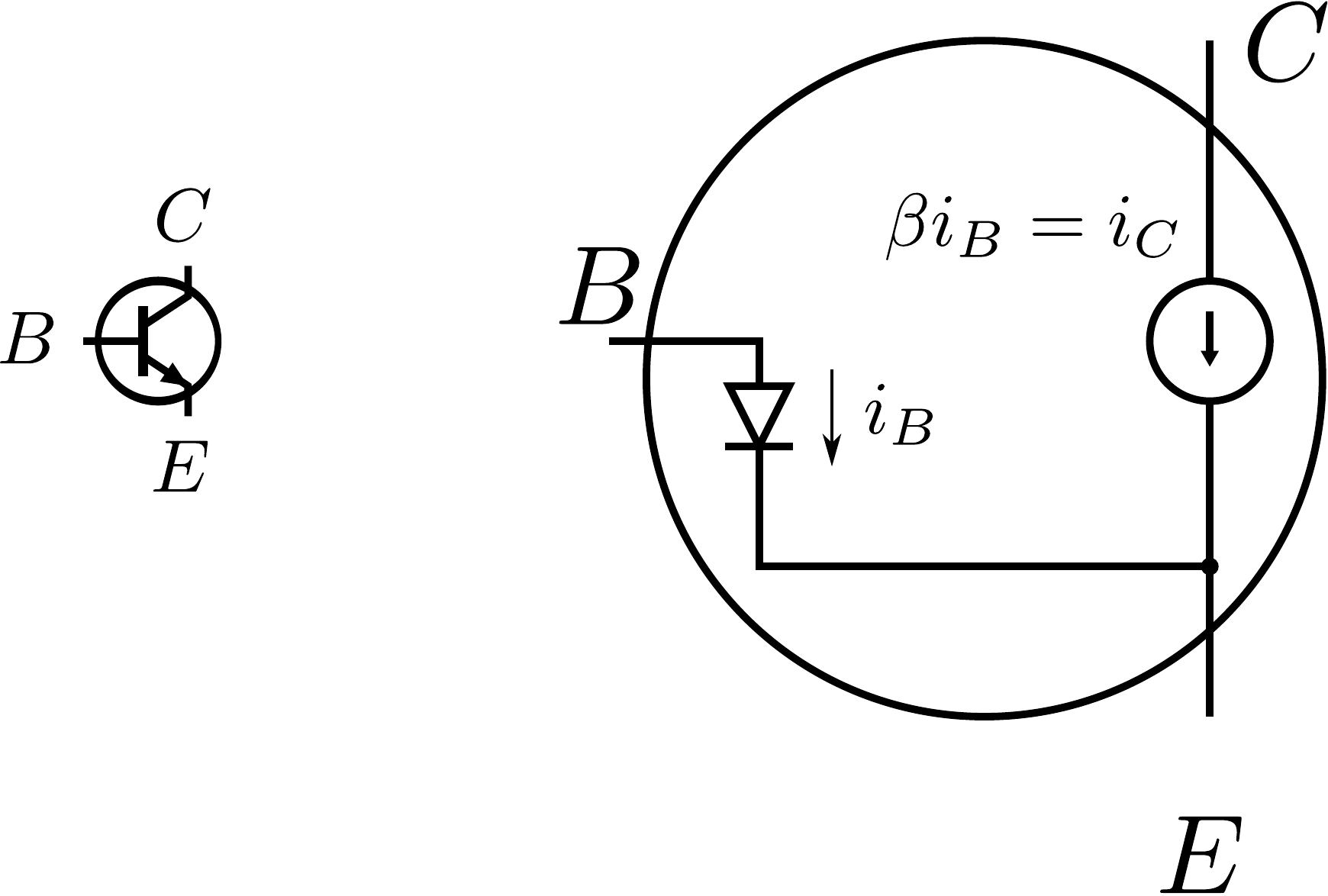}
\caption{NPN transistor model.}
\label{fig:simpleModel}
\end{figure}

In this section we derive an electronic implementation of the diagram in Fig.~\ref{FIG: burst circuit}.  The main active component is an NPN transistor. Its constitutive relations are given by Ebers-Moll equations~\cite{sedra} but, to simplify the analysis, we will model the transistor as a current source with current proportional to that of a diode standing at the base (see Fig.~\ref{fig:simpleModel}). Furthermore, we regard the diode as a perfect switch that opens whenever the base--emitter voltage is lower or equal to 0.6\:V, and closes otherwise. Also, the minimal voltage between the collector and the emitter is equal to 0.1\:V, a point at which the proportional relation among the base and emitter currents is lost.

\subsection{Non-smooth analysis}

The simplified transistor description leads to piecewise linear models that are easily dealt with using non-smooth analysis.
To make the paper self contained, we recall some definitions and results. See~\cite{makela,clarke} for further details.

\begin{definition}\cite[p. 27]{clarke}
 Consider a function $f: X \to \RE$ with $X$ a Banach space. The \emph{generalized directional derivative} of $f$ at $x$ in the direction of $\nu$, denoted
 $f^o(x;\nu)$, is defined as
 \begin{displaymath}
  f^o(x;\nu) = \limsup_{\substack{y \to x \\ t \downarrow 0}} \frac{f(y+t\nu) - f(y)}{t} \;.
 \end{displaymath}
 The \emph{generalized gradient} of $f$ at $x$, denoted $\partial f(x)$, is a subset of the dual
 space $X^*$ given by
 \begin{displaymath}
  \left\{ \xi \in X^* \: | \: f^o(x; \nu) \ge \langle \xi, \nu \rangle \; \text{for all $\nu$ in $X$} \right\} \;.
 \end{displaymath}
\end{definition}

Here we will set $X = \RE^n$ and identify $\RE^n$ with its dual, so $\partial f(x)$ is taken as a subset of $\RE^n$.
In this case we have:
\begin{theorem}\cite[p. 63]{clarke}
 Let $f$ be Lipschitz near $x$ and let $\Omega_f$ be the set of points at which $f$ fails to be differentiable.
 Suppose that $S$ is any set of Lebesgue measure 0 in $\RE^n$. Then,
 \begin{displaymath}
  \partial f(x) = \co \left\{ \lim \nabla f (x_i) \:|\: x_i \to x, \; x_i \notin S, \; x_i \notin \Omega_f \right\} \;.
 \end{displaymath}
\end{theorem}

\begin{theorem}\cite[p. 38]{makela}
 If $f$ is locally Lipschitz at $x$ and attains its extremum at $x$, then
 \begin{displaymath}
  0 \in \partial f(x) \;.
 \end{displaymath}
\end{theorem}

The smooth version of the implicit function theorem is a common tool for 
finding singular points in classical bifurcation theory. Let us present a non-smooth
version for non-smooth problems. Consider a Lipschitz function $H : \RE \times \RE \to \RE$
together with the equation
\begin{equation} \label{eq:implicit}
 H(y,z) = 0 \;.
\end{equation}
Suppose there is a pair 
$(\hat{y},\hat{z})$ that solves~\eqref{eq:implicit} and
\begin{displaymath}
 0 \notin \partial_z H(\hat{y},\hat{z}) \;, 
\end{displaymath}
where $\partial_z H$ is the generalized gradient with respect to $z$.
The implicit function theorem~\cite[p. 256]{clarke} states that
there exists a neighborhood $Y$ of $\hat{y}$ and a Lipschitz function $\zeta : \RE \to \RE$
such that $\hat{z} = \zeta(\hat{y})$ and such that, for every $y \in Y$,
\begin{displaymath}
 H(y,\zeta(y)) = 0 \;.
\end{displaymath}
\subsection{Voltage-controlled non-monotone characteristic} \label{sec:nm}

\begin{figure}
\centering
 \includegraphics[scale=0.29]{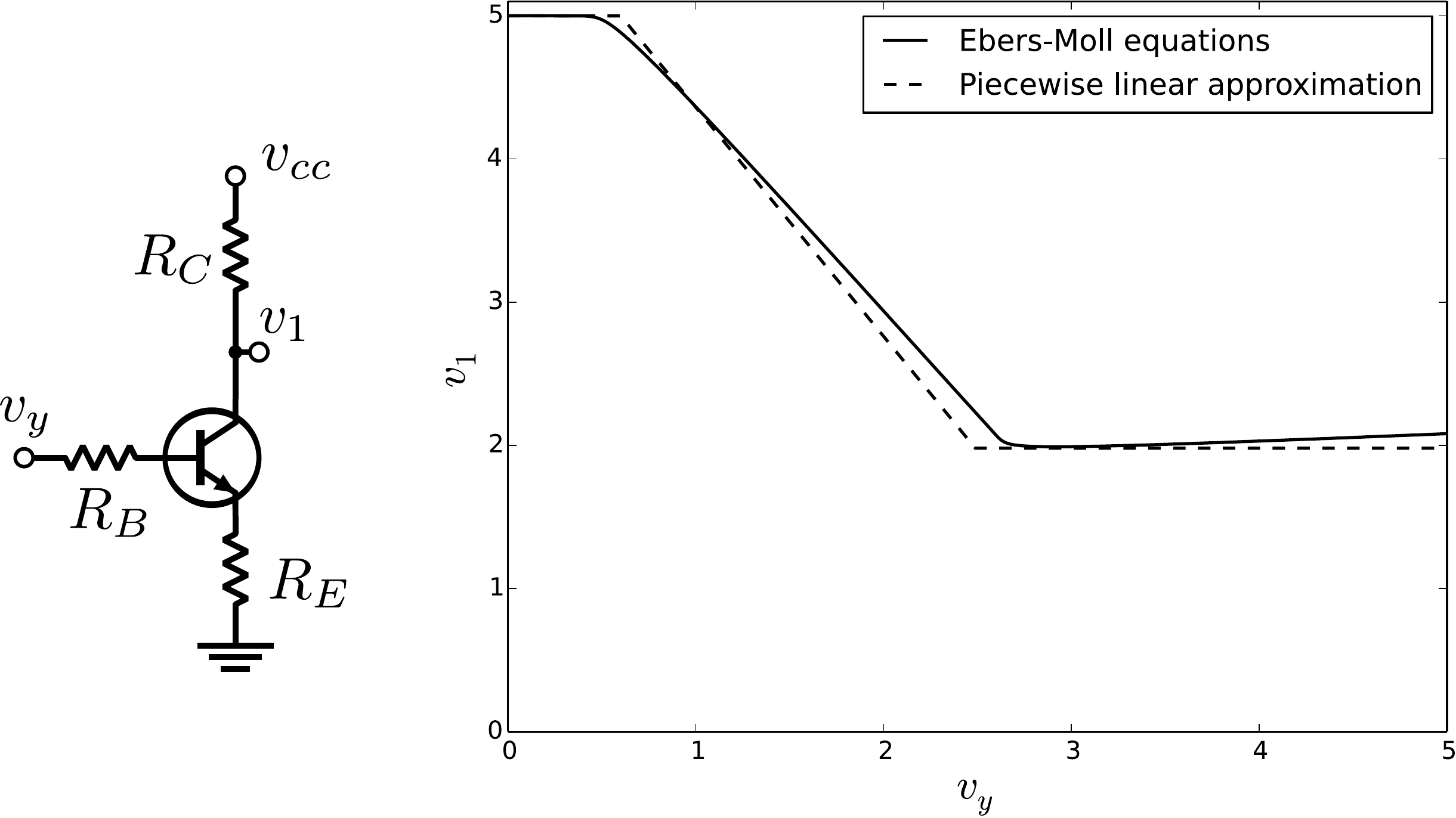}
\caption{Basic saturation. Common-emitter configuration (left). Piecewise linear approximation vs. exact solutions (right).}
\label{fig:ce_sat_plot}
\end{figure}

The non-monotone block of the circuit in Fig.~\ref{FIG: burst circuit} is realized as the difference of two monotone nonlinearities~\cite{Franci2014a}.
A common emitter NPN transistor configuration serves as a simple, natural saturation monotone nonlinearity (see Fig.~\ref{fig:ce_sat_plot}).
In view of the transistor model described above, the $v_y$--$v_1$ characteristic takes the piecewise linear form
\begin{equation}
\label{eq:sat}
 v_1 = v_{cc} - \Proj_S(g_1 (v_y-0.6)) \;,
\end{equation}
where $v_{cc} = 5\:V$ is the voltage at the power source,
\begin{displaymath}
 g_1 = \frac{\beta R_C}{R_B + (\beta + 1)R_E}
\end{displaymath}
is the voltage gain (the slope of the saturation), $\beta \approx 100$ is the transistor's current gain and
$S_1 = [0,v_{s}]$ with
\begin{displaymath}
 v_{s_1} = (v_{cc} - 0.1)\frac{R_C}{R_C + R_E}
\end{displaymath}
the saturation voltage. Usually, one chooses $R_B \ll \beta R_E$ so that $g_1 \approx R_C/R_E$, that is,
so that the dependence of $g_1$ on $\beta$ is negligible.

$\Proj_{S_1}$ is the operator that projects its argument into the set $S_1$, that is, 
\begin{displaymath}
 \Proj_{S_1}(v) = \argmin_{w \in S_1} \|v-w\| \;.
\end{displaymath}
For our particular $S_1$, the projection translates into
\begin{displaymath}
 \Proj_{S_1}(v) = 
  \begin{cases}
         0 & \text{if } v \le 0 \\
         v & \text{if } 0 \le v \le v_{s_1} \\
   v_{s_1} & \text{if } v_{s_1} \le v 
  \end{cases} \;.
\end{displaymath}
The projection captures the fact that the diode does not conduct when $v_{BE}$ is below 0.6\:V and that $i_C$ saturates
when $v_{CE}$ reaches 0.1\:V. To asses the quality of our estimation, we choose a set of parameters and compare~\eqref{eq:sat}
against the simulation obtained using {\tt ngspice} (which implements Ebers-Moll equations). The results are shown in Fig.~\ref{fig:ce_sat_plot}.

\begin{figure}
\centering
 \includegraphics[scale=0.29]{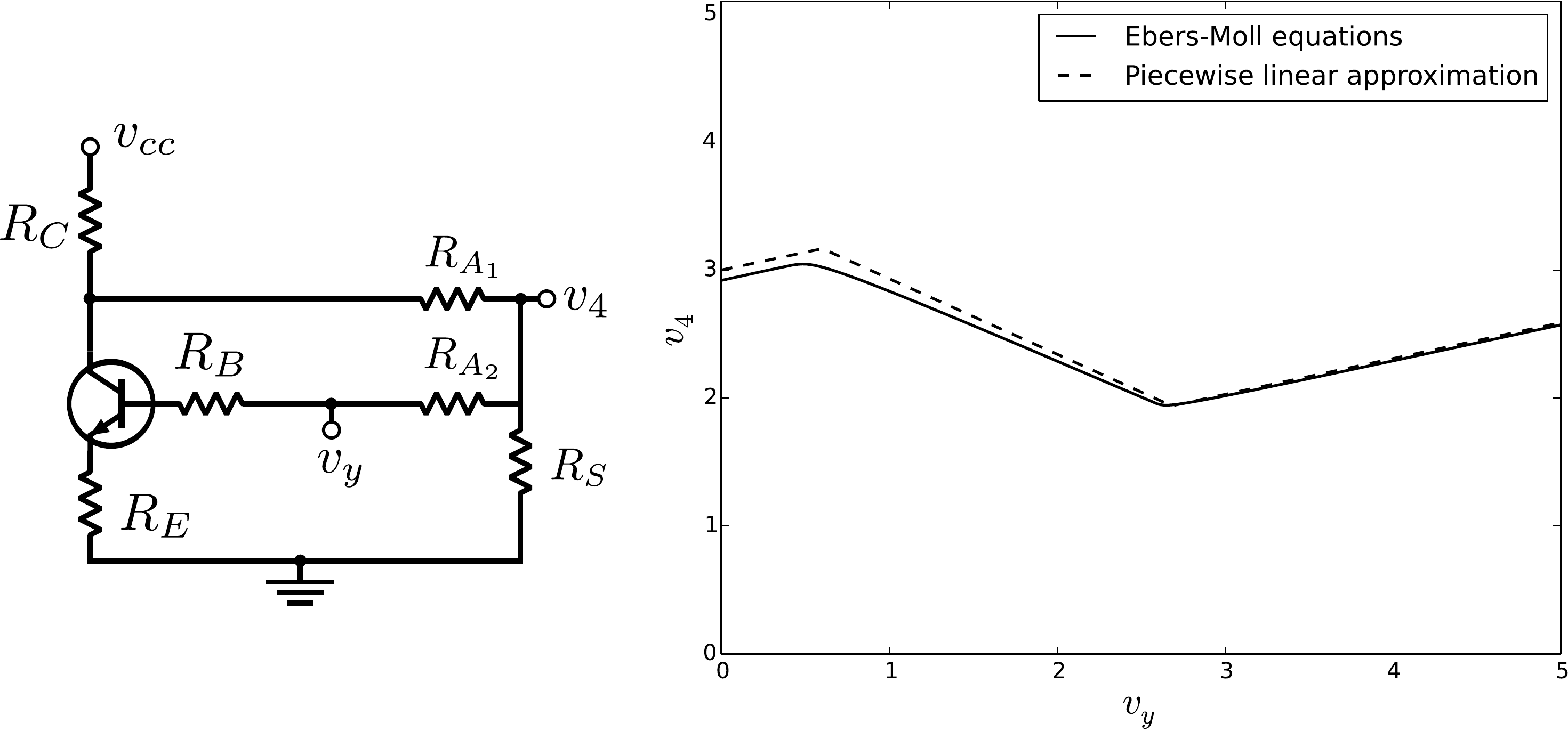}
\caption{Non monotone characteristic. Parallel interconnection of a saturation and
 a linear gain (left). Piecewise linear approximation vs. exact solutions (right).}
\label{fig:nm_nm_plot}
\end{figure}

The parallel interconnection achieving the non-monotone behavior is shown in Fig.~\ref{fig:nm_nm_plot}, left.
Applying Kirchhoff's laws and the piecewise linear model for the transistor one obtains
\begin{equation} \label{eq:nm}
 v_4 = g_2(v_{cc}-\Proj_S(g_1(v_y-0.6)))+g_3 v_y
\end{equation}
with
\begin{align*}
 g_2 &= \frac{R_S R_{A_2}}{R_S\left( R_{A_1} + R_{A_2} \right) + R_{A_1}R_{A_2}} \\
 g_3 &= \frac{R_S R_{A_1}}{R_S\left( R_{A_1} + R_{A_2} \right) + R_{A_1}R_{A_2}}
\end{align*}
(see Fig.~\ref{fig:nm_nm_plot}, right).

The characteristic~\eqref{eq:nm} is non-monotone whenever strict extrema are present. Recall that a
necessary condition for the presence of extrema is $0 \in \partial v_4$,
where $\partial v_4$ is the generalized gradient of $v_4$ with respect to $v_y$.
Note that
\begin{equation} \label{eq:Psi}
 \partial \Proj_{S_1}(v) = 
  \Psi_{S_1}(v) :=
  \begin{cases}
       \left\{ 0 \right\} & \text{if } v \notin S_1 \\
                    [0,1] & \text{if } v \in \partial S_1\\
       \left\{ 1 \right\} & \text{if } v \in \Int S_1
  \end{cases} 
\end{equation}
with $\Int S_1$ and $\partial S_1$ the interior and the boundary of $S_1$, respectively.
Application of the chain rule to~\eqref{eq:nm} gives
\begin{displaymath}
 \partial v_4 = g_3 - g_1g_2\Psi_{S_1}(g_1(v_y-0.6)) \;.
\end{displaymath}
A necessary condition for $0 \in \partial v_4$ is 
\begin{equation} \label{eq:g1g2}
 g_1 g_2 \ge g_3 \;.
\end{equation}
This inequality imposes a set of conditions to be satisfied by the resistors
$R_{A_1}$, $R_{A_2}$, $R_{S}$, $R_B$, $R_C$ and $R_E$.

When inequality~\eqref{eq:g1g2} holds strictly (the equality does not happen in practice), the inclusion
$0 \in \partial v_4$ occurs at the boundary of $S_1$, that is, at the points $v_y$ that satisfy
\begin{displaymath}
 g_1(v_y - 0.6) = 0 \quad \text{and} \quad g_1(v_y - 0.6) = v_{s_1} \;.
\end{displaymath}
We thus have the following points of interest for~\eqref{eq:nm}:
\begin{align*}
                                    v_4(0.6) &= g_2 v_{cc} + g_3 0.6 \quad (\text{a local maximum}) \\
 v_4\left( \frac{v_{s_1}}{g_1} + 0.6 \right) &= g_2 (v_{cc} - v_{s_1}) \\
                                              &{} \: + g_3 \left( \frac{v_{s_1}}{g_1} + 0.6 \right) \quad (\text{a local minimum}) \\
                                      v_4(0) &= g_2 v_{cc} \quad (\text{initial point}) \\
                                 v_4(v_{cc}) &= (g_2 + g_3) v_{cc} - g_2 v_{s_1} \quad (\text{final point}) \;.
\end{align*}
The computation of such points is again useful for choosing the appropriate resistors. The objective is to have the minimum occur 
as close as possible to $v_y = v_{cc}/2 = 2.5\:V$, and to have the largest possible excursion along the $v_4$ axis. 

\begin{figure}
\centering
 \includegraphics[scale=0.29]{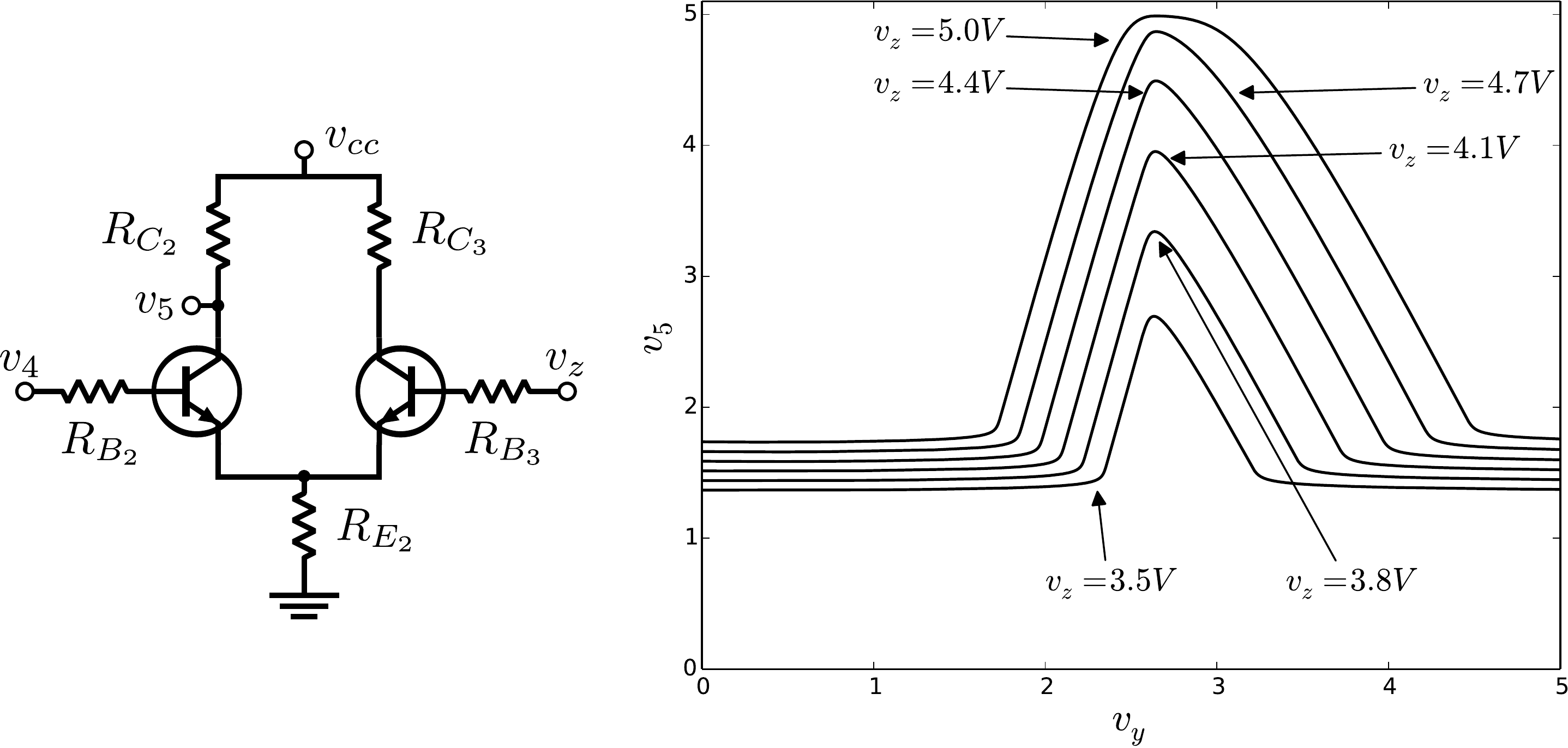}
\caption{The voltage-controlled non-monotone characteristic is realized by cascading a differential amplifier (left)
 with a fixed non-monotone characteristic (see Fig.~\ref{fig:nm_nm_plot}). Input--output response for different values of $v_z$ (right).}
\label{fig:diff_mod_nm_plot}
\end{figure}

In the block realization in Fig.~\ref{FIG: burst circuit}, the output of the non-monotone block is modulated by the unfolding parameter $\alpha$.
Such modulation is necessary to transition along the three possible
characteristics in Fig.~\ref{FIG:bursting phase portrait}-C.
The same modulation is achieved here by cascading a differential amplifier with the non-monotone characteristic.
In this way, an external voltage can be used to the scale the input--output characteristic, as shown in
Fig.~\ref{fig:diff_mod_nm_plot}, left. We have
\begin{align*}
 v_4 - 0.6 &= \left( R_{B_2} + \bar{R}_{E} \right)i_{B_2} + \bar{R}_{E} i_{B_3} \\
 v_z - 0.6 &= \bar{R}_{E} i_{B_2} + \left( R_{B_3} + \bar{R}_{E} \right) i_{B_3}
\end{align*}
with $\bar{R}_{E} = (\beta+1)R_{E}$. The currents are thus
\begin{displaymath}
 \begin{pmatrix}
  i_{B_2} \\ i_{B_3}
 \end{pmatrix}
 =
 \frac{1}{d}
 \begin{pmatrix}
   \bar{R}_{E} + R_{B_3} & - \bar{R}_{E} \\
  -\bar{R}_{E}           &  \bar{R}_{E} + R_{B_2} 
 \end{pmatrix}
 \begin{pmatrix}
  v_4 - 0.6 \\ v_z - 0.6 
 \end{pmatrix}
\end{displaymath}
with
\begin{displaymath}
 d = \bar{R}_{E}\left( R_{B_2} + R_{B_3} \right) + R_{B_2}R_{B_3} \;.
\end{displaymath}
Let $\bar{R}_{C_2} = \beta R_{C_2}$ and $\bar{R}_{C_3} = \beta R_{C_3}$. When $i_{C_3}$
does not saturate, i.e., when
\begin{displaymath}
 \bar{R}_E i_{B_2} + (\bar{R}_{C_3} + \bar{R}_E) i_{B_3} < v_{cc} - 0.1 \;,
\end{displaymath}
the output voltage is
\begin{equation} \label{eq:diff}
 v_5 = v_{cc} - \Proj_{S_2}\left(g_4 (v_4-0.6) - g_5 (v_z-0.6) \right) \;,
\end{equation}
where
\begin{align*}
 g_4 &= \frac{1}{d} \bar{R}_{C_2}\left( \bar{R}_{E} + R_{B_3} \right) \\
 g_5 &= \frac{1}{d} \bar{R}_{C_2} \bar{R}_{E} \;.
\end{align*} 
and the interval $S_2 = [0,v_{s_2}]$ is determined by
\begin{displaymath}
 v_{s_2} = \frac{\bar{R}_{C_2}}{\bar{R}_E + \bar{R}_{C_2}}\left( v_{cc}-0.1-\bar{R}_E i_{B_3} \right) \;.
\end{displaymath}

The voltage $v_4$ is called the non-inverting input
and $v_z$ the inverting one. The complete block is shown in Fig.~\ref{fig:mod_mirror} and the input--output
characteristic is shown in Fig.~\ref{fig:diff_mod_nm_plot}, right.

\subsection{Hysteretic characteristic}

\begin{figure}
\centering
 \includegraphics[scale=0.29]{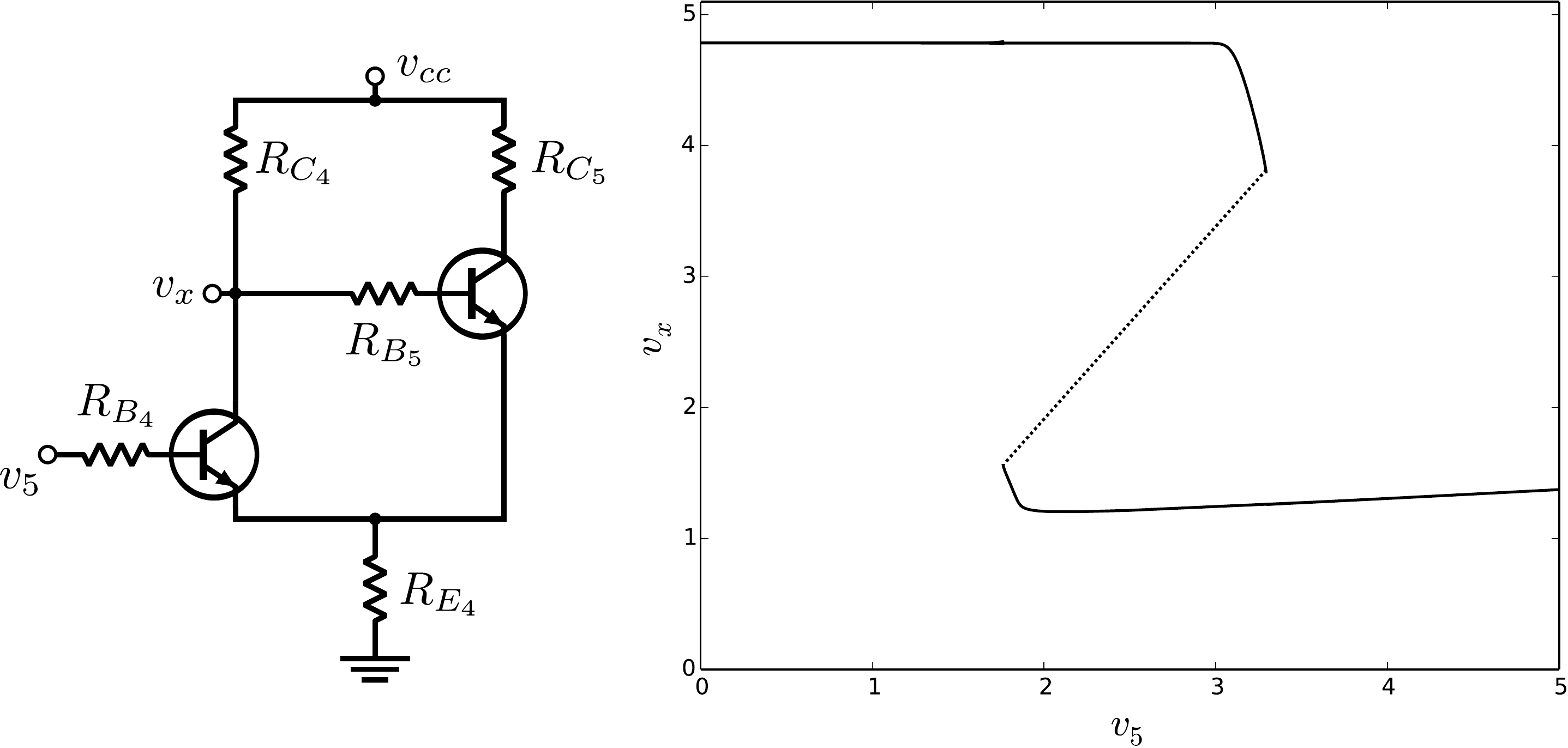}
\caption{Hysteretic characteristic. Input--output response. To illustrate the hysteretic behavior, the input
 $v_5$ is swept from 0 to 5V and back. A dotted line has been manually added to sketch the solutions that
 the numerical solver cannot find and which correspond to unstable steady states.}
\label{fig:h_cir_plot}
\end{figure}

The hysteretic block is built as the positive feedback of a basic saturation and a linear gain~\cite{Franci2014a}.
This is achieved at once with another differential amplifier. By letting $v_5$ be the inverting input, $v_6$ the 
non-inverting input and $v_x$ the output, we obtain, cf.~\eqref{eq:diff},
\begin{displaymath}
 v_x =  v_{cc} - \Proj_{S_3}\left(g_6 (v_5-0.6) - g_7 (v_6-0.6) \right) \;.
\end{displaymath}
Positive feedback is then achieved simply by setting $v_6 = v_x$, as shown in Fig.~\ref{fig:h_cir_plot}. This results in
the piecewise linear characteristic
\begin{multline} \label{eq:h}
 F(v_5,v_x) = v_x - v_{cc} \\ + \Proj_{S_3}\left( g_6(v_5-0.6) - g_7(v_x-0.6) \right) = 0 \;.
\end{multline}
It follows from the implicit function theorem that a necessary condition for the existence of
singular points is $0 \in \partial F(v_5,v_x)$, with the generalized gradient taken with
respect to $v_x$. Application of the chain rule to~\eqref{eq:h} gives
\begin{displaymath}
 \partial F(v_5,v_x) = 1 - g_7 \Psi_{S_3}\left( g_6(v_5-0.6) - g_7(v_x-0.6) \right) \;.
\end{displaymath}
A necessary condition for $0 \in \partial F(v_5,v_x)$ is thus $g_7 \ge 1$. There are two points of
singularity occurring at the boundary of $S_3$. The first one is characterized by
\begin{displaymath}
 g_6(v_5-0.6) = g_7(v_x-0.6)
\end{displaymath}
which, together with the condition $F(v_5,v_x) = 0$, gives
\begin{displaymath}
 v_x = v_{cc} \quad \text{and} \quad v_5 = \frac{g_6}{g_7}v_{cc} + \frac{g_6 - g_7}{g_6}0.6 \;.
\end{displaymath}
The other point of singularity is determined by 
\begin{displaymath}
 g_6(v_5-0.6) = g_7(v_x-0.6) + v_{s_3} \;,
\end{displaymath}
which gives
\begin{displaymath}
 v_x = v_{cc}-v_{s_3} \quad \text{and} \quad v_5 = \frac{g_7}{g_6}v_{cc} + \frac{g_6 - g_7}{g_6}0.6 - \frac{g_7 - 1}{g_6}v_{s_3}
\end{displaymath}
(see Fig.~\ref{fig:h_cir_plot}). For this block, the resistors were chosen so that the first and second singularities occur,
respectively, at one and two thirds of the chosen voltage range of 5\:V.

\subsection{Voltage-controlled mirrored hysteresis}

\begin{figure*}
\centering
 \includegraphics[scale=0.29]{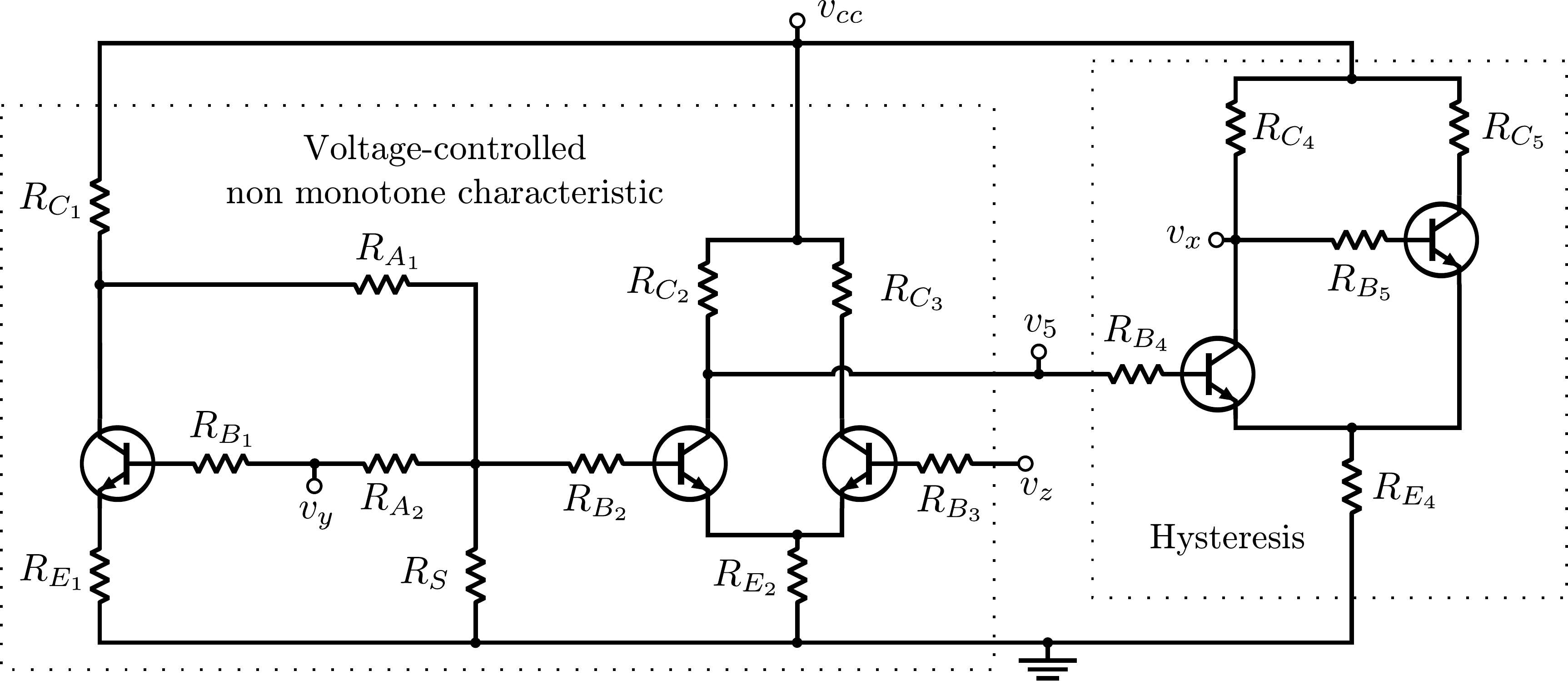}
\caption{Circuit realizing the mirrored hysteresis in simulation. By suitably changing $v_z$, the $v_y$--$v_x$ characteristic
can adopt the three forms portrayed on Fig.~\ref{fig:mod_mirror_plot} (cf.  the right part of Fig.~\ref{FIG:bursting phase portrait}).}
\label{fig:mod_mirror}
\end{figure*}

\begin{figure*}
\centering
 \begin{subfigure}[t]{0.32\textwidth}
  \centering
  \includegraphics[width=0.95\columnwidth]{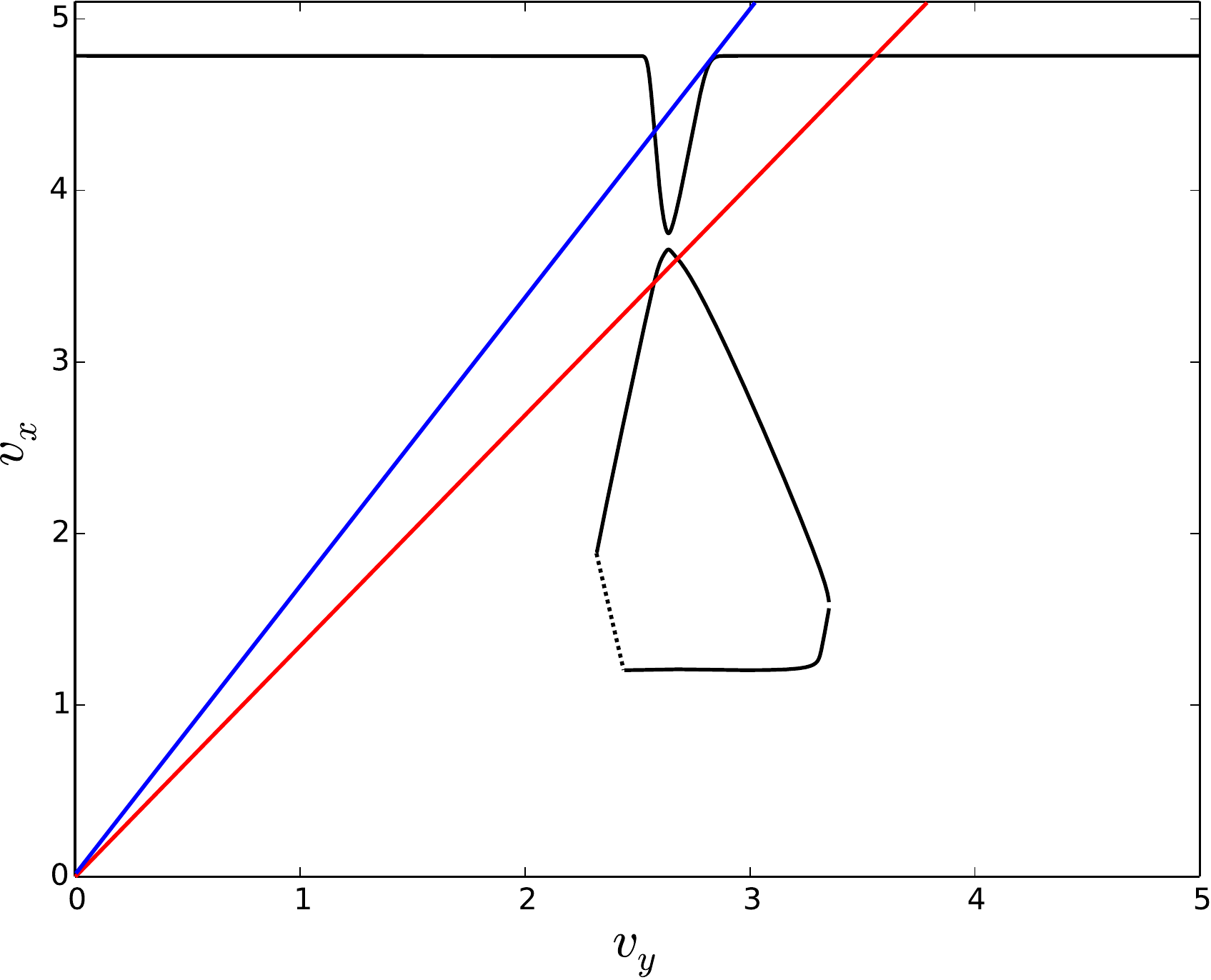}
  \caption{$v_z = 3.8\:V$. A saddle point, an unstable and a stable node are present (red). Only a stable
   node is present (blue). In both cases, almost all trajectories converge to the stable node.}
  \label{fig:mod_mirror_1_plot}
 \end{subfigure} \hspace{2pt}
 \begin{subfigure}[t]{0.32\textwidth}
  \centering
  \includegraphics[width=0.95\columnwidth]{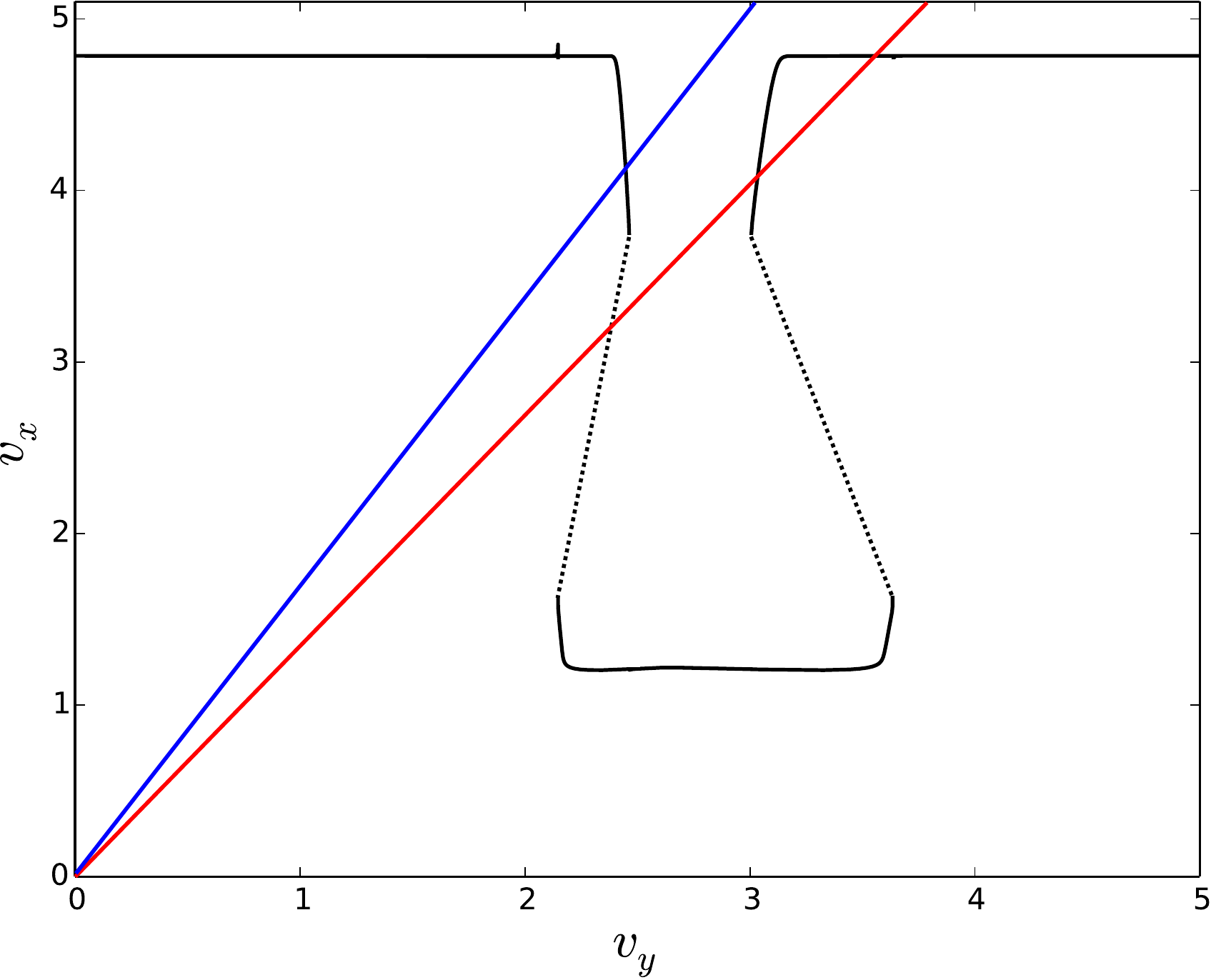}
  \caption{$v_z = 4.1\:V$. A saddle point, a stable limit cycle around an unstable node and a stable node coexist (red).
   The only attractor is a stable node (blue).}
  \label{fig:mod_mirror_2_plot}
 \end{subfigure} \hspace{2pt}
 \begin{subfigure}[t]{0.32\textwidth}
  \centering
  \includegraphics[width=0.95\columnwidth]{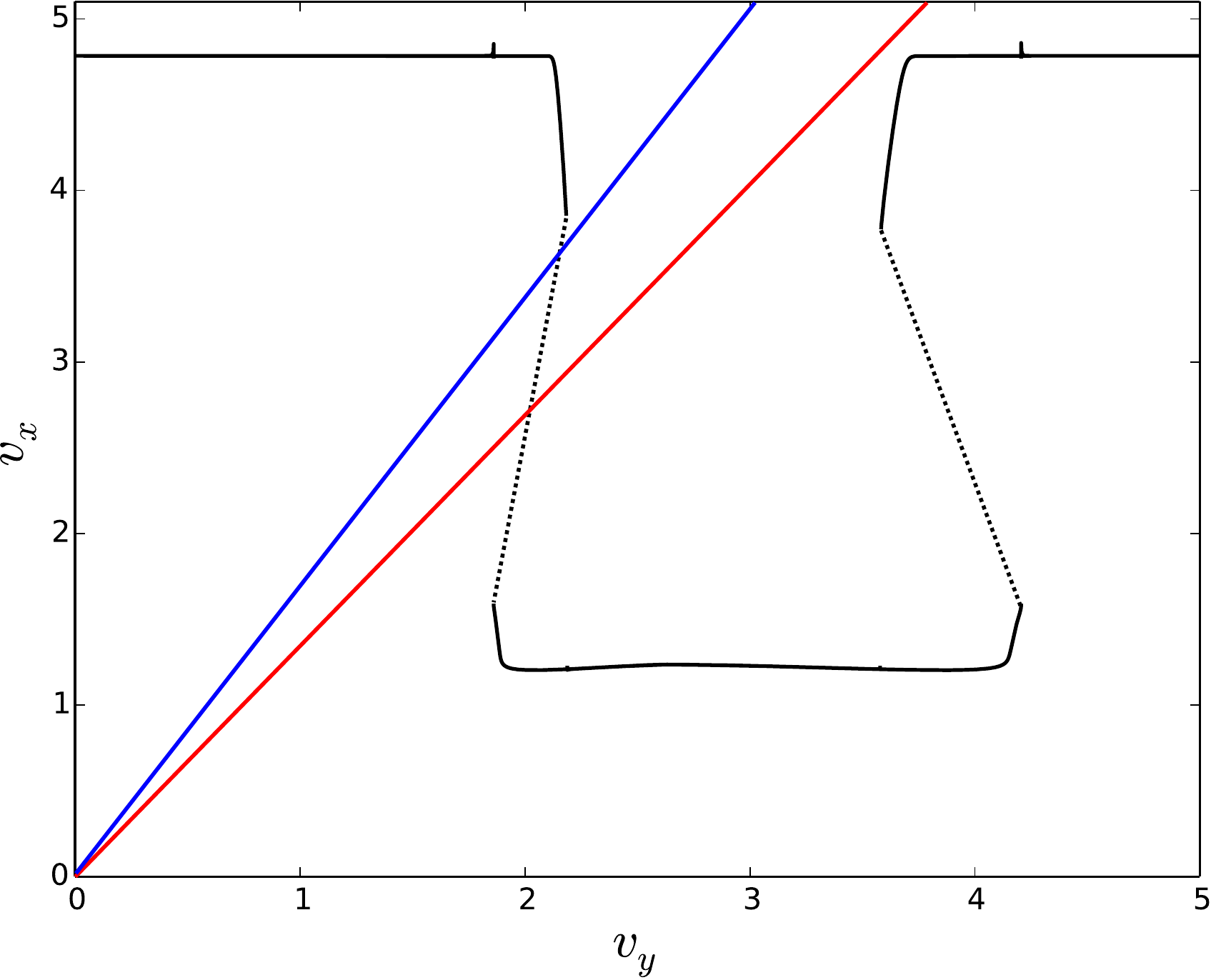}
  \caption{$v_z = 4.7\:V$. A stable limit cycle around an unstable node exists. Almost all trajectories converge to
   the limit cycle (red and blue).}
  \label{fig:mod_mirror_3_plot}
 \end{subfigure}
\caption{Mirrored-hysteresis. Solid black lines correspond to the $x$-steady states found by solving the circuit with {\tt ngspice}.
 Dotted lines are manually added to sketch the solutions not found by the solver.
 Red and blue lines correspond to the desired $v_y$-nullclines.}
\label{fig:mod_mirror_plot}
\end{figure*}

The cascade of the voltage-controlled non-monotone block and the hysteresis is shown in Fig.~\ref{fig:mod_mirror}. This interconnection
realizes the block diagram in Fig.~\ref{FIG: wcusp circuit}. The circuit establishes the desired static behavior relating
$v_x$, $v_y$ and $v_z$, that is, it implements a voltage characteristic that is topologically equivalent to the algebraic variety 
given by~\eqref{EQ:bursting normal form critical manifold}. In other words, it produces the desired voltage-controlled mirror
hysteresis (see Figs.~\ref{fig:mod_mirror_plot}, black traces).
Compared to Figs.~\ref{FIG:bursting phase portrait}-C, the simulated responses are, on one hand, stiffer and, on the other, reflected along
a horizontal axis. However, this does not alter the qualitative picture in terms of number and type of different possible attractors. In
this sense we say that the two portraits are qualitatively equivalent.

\subsection{Burster}

We now transform the static circuit in Fig.~\ref{fig:mod_mirror} into a dynamic circuit exhibiting the same qualitative
dynamics as model~\eqref{EQ:bursting normal form}.
The parasitic capacitances of the transistors provide the fast $v_x$ dynamics and set its corresponding timescale (cf.
$H_f$ in Fig.~\ref{FIG: burst circuit}). The voltage $v_x$ is fed back to $v_y$ by means of a resistive voltage divider
and a capacitor. The values of the resistors and the capacitor determine the timescale of the slow $v_y$ dynamics as well
as the slope of its nullcline (cf. $H_s$ in Fig.~\ref{FIG: burst circuit}). For bursting, we choose the resistors in such a way
that this slope is small and the $v_y$ nullcline intersects both branches of the mirrored hysteresis, as in Fig.~\ref{FIG:bursting phase portrait}.

\begin{figure}
\centering
 \includegraphics[width=0.90\columnwidth]{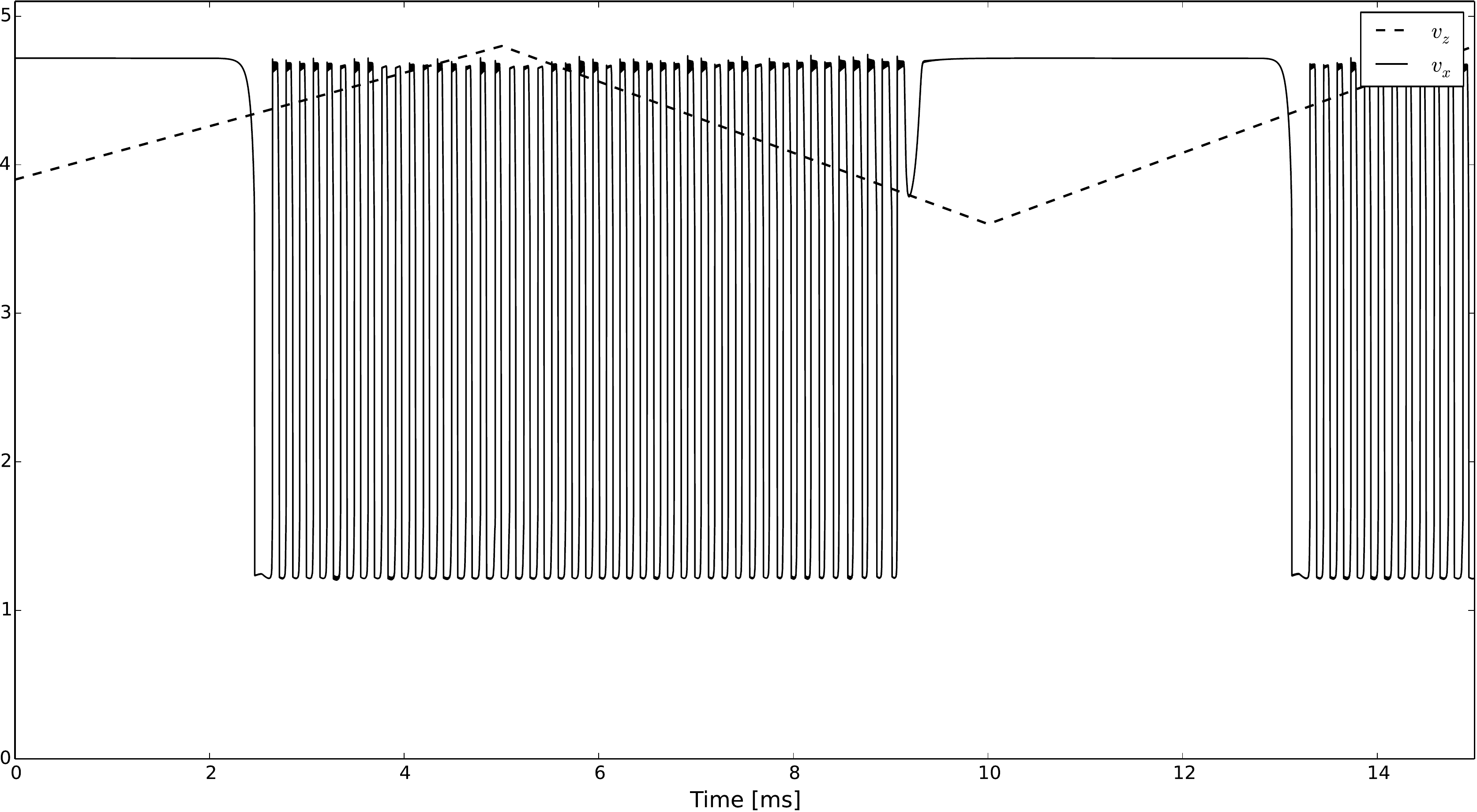}
\caption{Transition between a stable node (constant output) and a stable limit cycle (oscillations) in the circuit of Fig.~\ref{fig:vyz}.
 The transitions occur at different values of $v_z$, which indicates that the stable node and the stable limit cycle coexist for some values of $v_z$.}
\label{fig:vy_plot}
\end{figure}

Fig.~\ref{fig:mod_mirror_plot} confirms the qualitative equivalence of the circuit in Fig.~\ref{fig:mod_mirror} and model~\eqref{EQ:bursting normal form}. By sweeping $v_z$ we recover the same qualitative phase portraits of Fig.~\ref{FIG:bursting phase portrait}-C left, which underlie the behavior simulated in Fig.~\ref{fig:vy_plot}.
A key ingredient in the bursting behavior is the bistability of the limit cycle and the node (Fig.~\ref{fig:mod_mirror_2_plot}, red). The presence of this phenomenon can be asserted by noting that  the transition from the constant output (the stable node) to the oscillating behavior (the limit cycle) occurs at a
higher value of $v_z$ than the one for the transition from the oscillating behavior back to the constant output.

\begin{figure}
\centering
 \includegraphics[width=0.90\columnwidth]{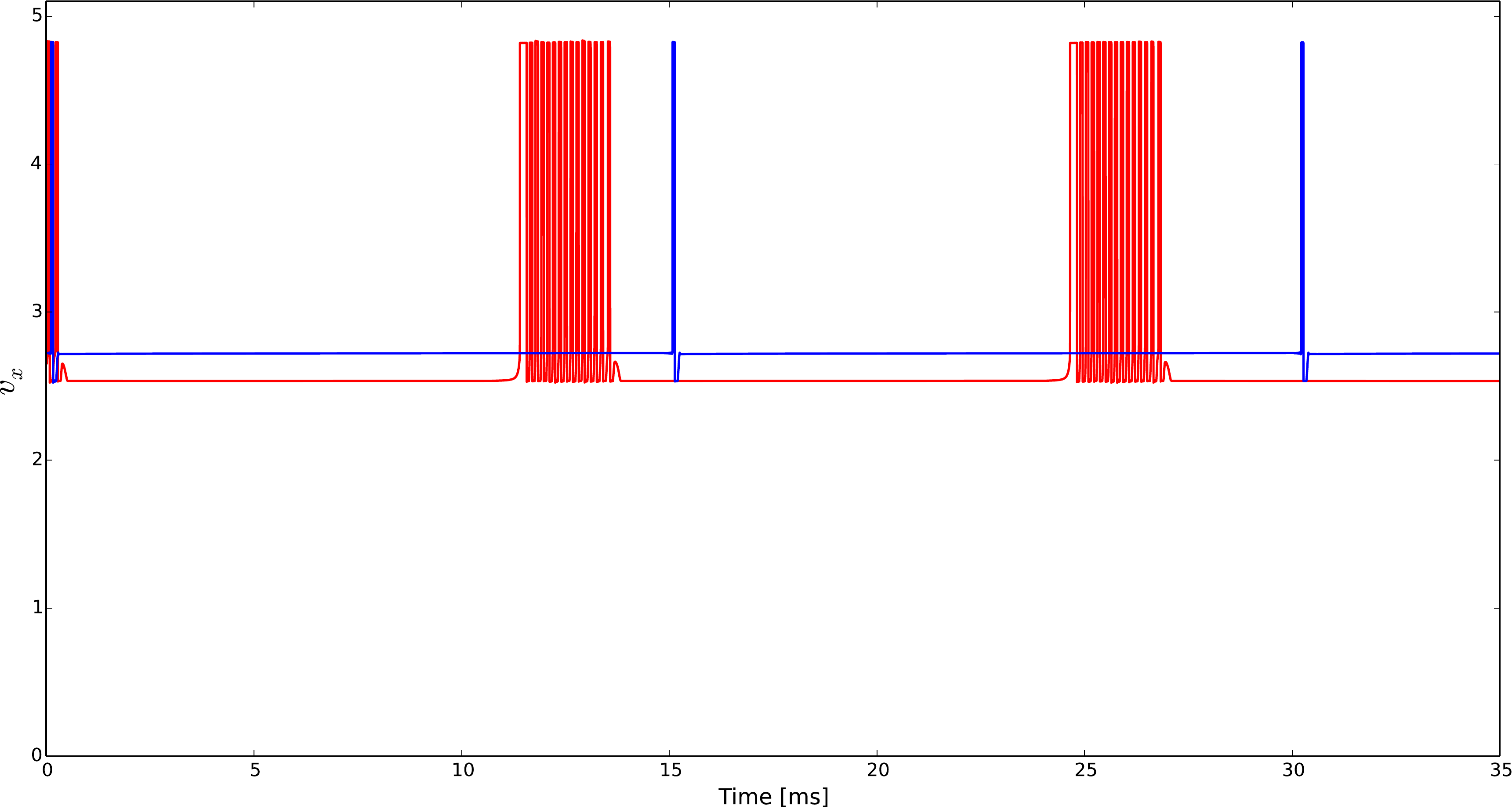}
\caption{Bursting (red) and tonic spiking (blue). Simulation.}
\label{fig:bs_plot}
\end{figure}

Bursting is finally achieved by feeding $v_x$ back to $v_z$ through the ultra slow filter $H_u$. This is again realized
with a voltage divider and a capacitor, but now the circuit's time-constant is chosen much larger. To ensure a robust operation,
the output of $H_u$ is amplified so that $v_z$ exhibits a large swing. In fact it is the complement of $v_x$ that is passed
through an amplifier with negative slope (this accounts for the sixth transistor). The complete circuit is shown in Fig.~\ref{fig:vyz}
and its time response is shown in Fig.~\ref{fig:bs_plot}.

\subsection{Spiker}

Recall that the mirrored hysteresis captures both modes of operation: bursting and tonic spiking. {Geometrically, the difference between the two behaviors is the locus of the stable fixed point, as sketched in Fig.~\ref{FIG:bursting phase portrait}. In our circuit, we recover the same geometric picture by changing
the slope of the $v_y$-nullcline via tuning of the resistance $R_{i_1}$ and $R_{i_2}$. When the slope of the $v_y$ nullcline is sufficiently large this line solely intersects the left branch of the mirrored hysteresis at the transition between spiking and resting
(see Figs.~\ref{fig:mod_mirror_2_plot} and~\ref{fig:mod_mirror_3_plot}, blue), thus destroying the possibility of the bistability underlying bursting. The model is in this case in the tonic spiking mode shown in Fig.~\ref{fig:bs_plot}.}

\begin{figure*}
\centering
 \includegraphics[scale=0.29]{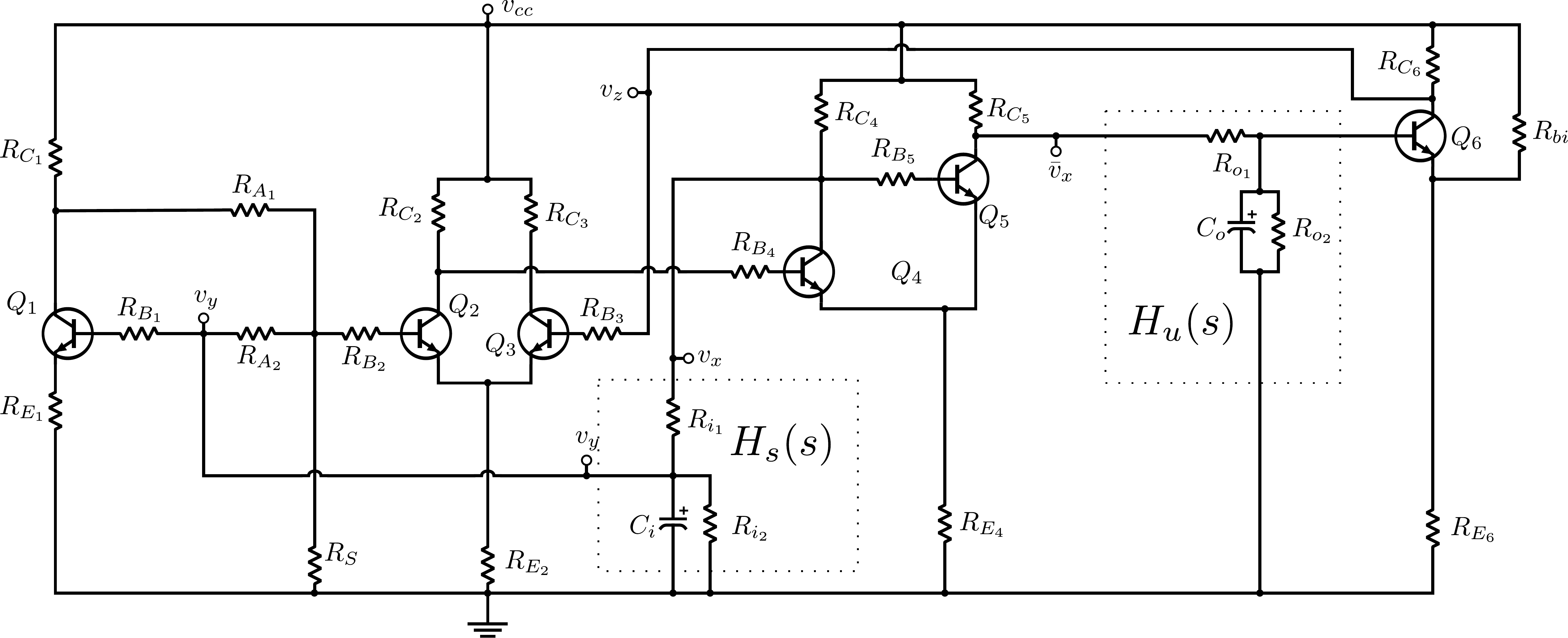}
\caption{Complete circuit capable of bursting and spiking in simulation. The ultra slow dynamics results from the capacitor $C_o$ and the voltage
 divider $R_{o_1}$, $R_{o_2}$. The slow dynamics results from the capacitor $C_i$ and the voltage divider $R_{i_1}$, $R_{i_2}$. 
 The fast dynamics results from the transistor's parasitic capacitances.}
\label{fig:vyz}
\end{figure*}

All the results of this section can be easily reproduced using the code provided in the appendix.

{
\subsection{Modulation of excitability properties in electronic devices}

The possibility of reliably switching between bursting and tonic spiking is relevant because it provides a means
to modulate the input--output behavior of our circuit. The output is the fast capacitor voltage $v_x$. We introduce
inputs by adding a current source $i_{\mathrm{app}}$ at the node labeled $v_z$ in Fig.~\ref{fig:vyz}. The injected current corresponds
to $u$ in~\eqref{EQ:bursting normal form}. Sufficiently large positive applied
currents set the circuit in a stable resting state, both in the tonic and bursting modes. The resting state is however
excitable: the response to excitatory (that is, negative for our circuit) input 
is large and highly nonlinear (spiking), reflecting the latent nonlinear dynamics of the circuit (Fig.~\ref{FIG:exc}).

\begin{figure}
\centering
\includegraphics[width=0.75\columnwidth]{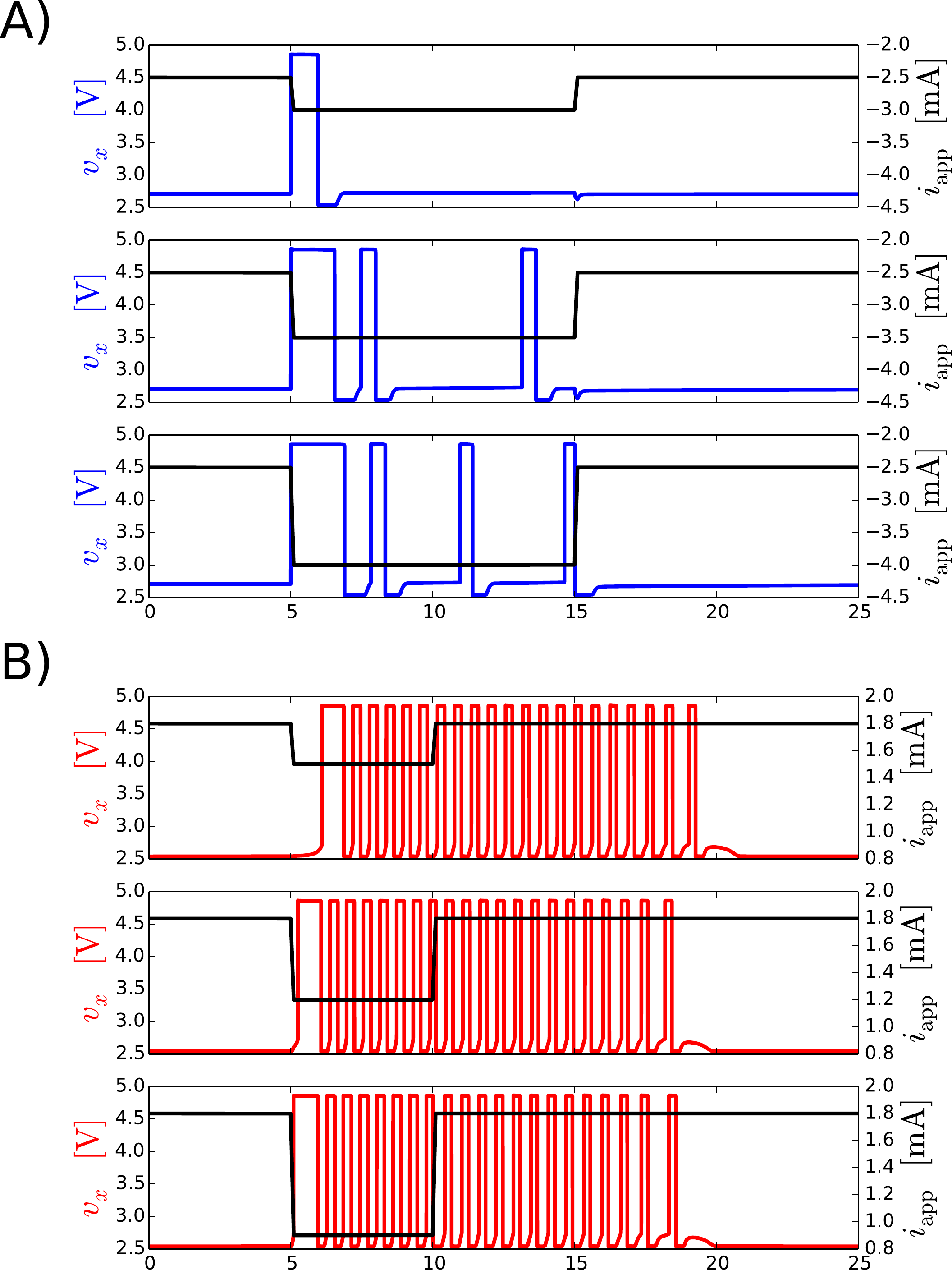}
\caption{Input--output response of the neuromorphic circuit in tonic (A) and bursting (B) modes. The input 
is provided by a source that injects a current $i_{\mathrm{app}}$ at the node labeled $v_z$ (Fig.~\ref{fig:vyz}) and is depicted in black. 
The response of the output voltage $v_x$ to different current steps when in tonic and bursting mode is depicted in blue
and red, respectively. }\label{FIG:exc}
\end{figure}

Excitability is sharply different in tonic and bursting modes \cite{Sherman2001}. The tonic mode is characterized by a quasi-linear coding of incoming inputs: the response of the circuit lasts only as long as the excitatory input is applied and the elicited spiking frequency is roughly proportional to the magnitude of the excitatory inputs. On the contrary, the bursting mode is characterized by a nonlinear detection mechanism of incoming inputs: the response of the circuit lasts for a fixed amount of time that can outlast the length of the excitatory stimuli (memory) and the elicited frequency is roughly independent of the input magnitude. Excitability in bursting mode serves as a bell ring signaling the arrival of new incoming inputs (for instance an unexpected sensory stimulus), whereas excitability in tonic mode serves as a frequency-coding mechanism to transmit information about those inputs. The possibility of switching between the two modes is widespread in the brain \cite{Sherman2001,lee2012neuromodulation}.

Our circuit exhibits the same behavioral transition. In the tonic mode the voltage output response roughly codifies the input current step magnitude (Fig.~\ref{FIG:exc}A). In the bursting mode it responds with bursts whose length and interspike frequency are virtually independent of the magnitude of the input current (Fig.~\ref{FIG:exc}B).

}
\section{Real-world implementation} \label{SEC: experimental}

The real circuit was built using the same methodology used for designing the simulated circuit.
First, a non-monotone block was built as in Section~\ref{sec:nm}. Using the formulas in that section,
the resistors were tuned so that the local minimum occurs at half the voltage range, $2.5\:V$. The
output of the non-monotone block was connected to a differential amplifier without much difference
with the simulated design. A hysteretic characteristic was implemented by means of a differential 
amplifier in positive feedback. Again, the resistors were chosen in order to have have a symmetric
characteristic with a large swing in $v_x$. The direct interconnection of the voltage-controlled block
and the hysteresis did not work as expected, the main reason being the fact that the hysteretical 
block draws a non-negligible current from the non-monotone block and changes its behavior (a phenomenon usually referred
to as \emph{loading}). To overcome this issue, an extra transistor in an emitter-follower configuration
was inserted between these two blocks. In this configuration, the extra transistor presents high input and
low output impedances, allowing the hysteretical block to be driven by the non-monotone block, but preventing
the hysteretical block from affecting the behavior of the non-monotone block. This effectively overcomes the
loading problem and simplifies the tuning of the circuit parameters, but induces a stiffer characteristic. 

\begin{figure*}
\centering
 \includegraphics[scale=0.35]{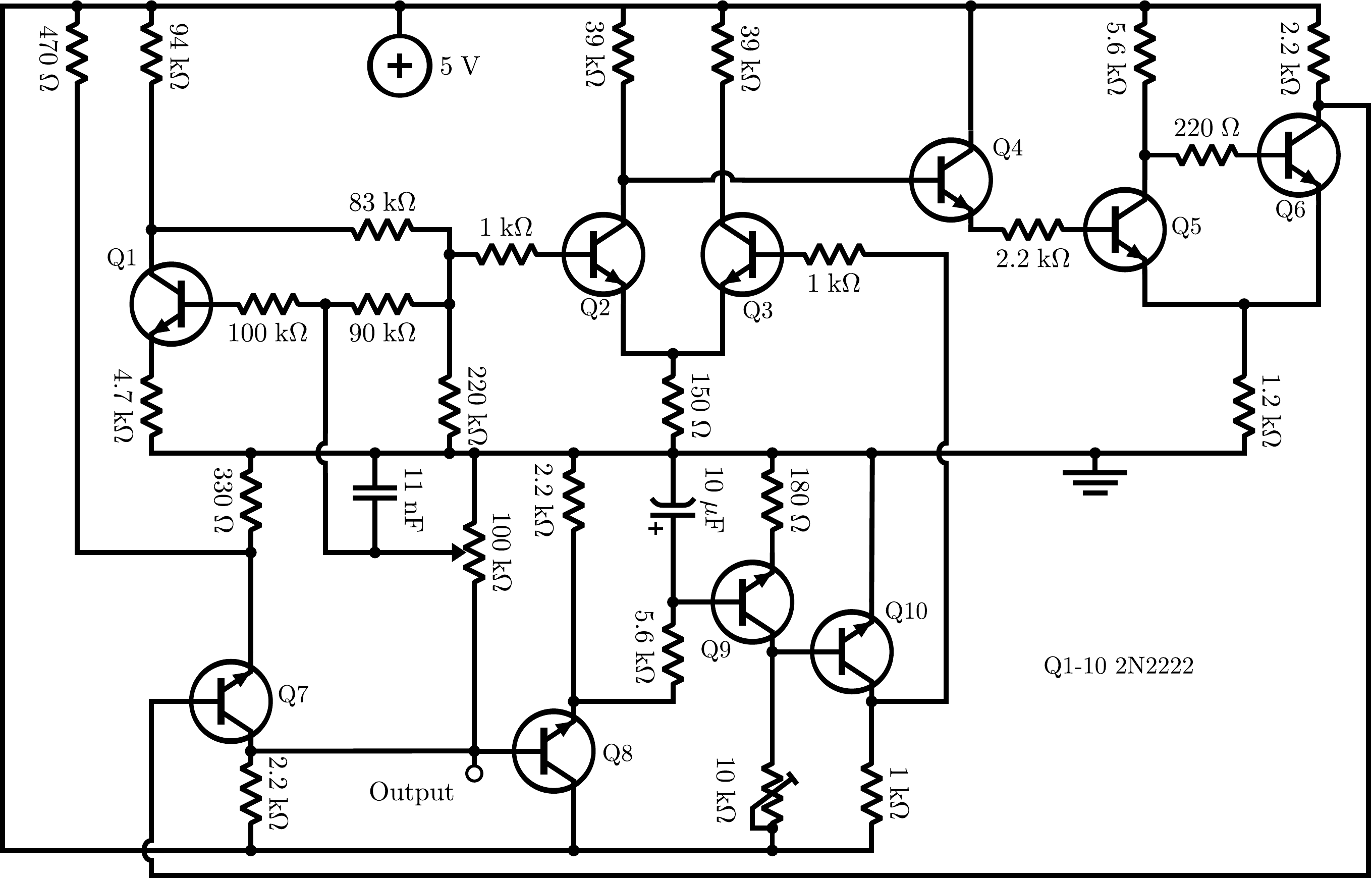}
\caption{Circuit used for experimental testing.}
\label{fig:final_circuit}
\end{figure*}

\begin{figure*}
\centering
 \begin{subfigure}[t]{0.32\textwidth}
  \centering
  \includegraphics[scale=0.30]{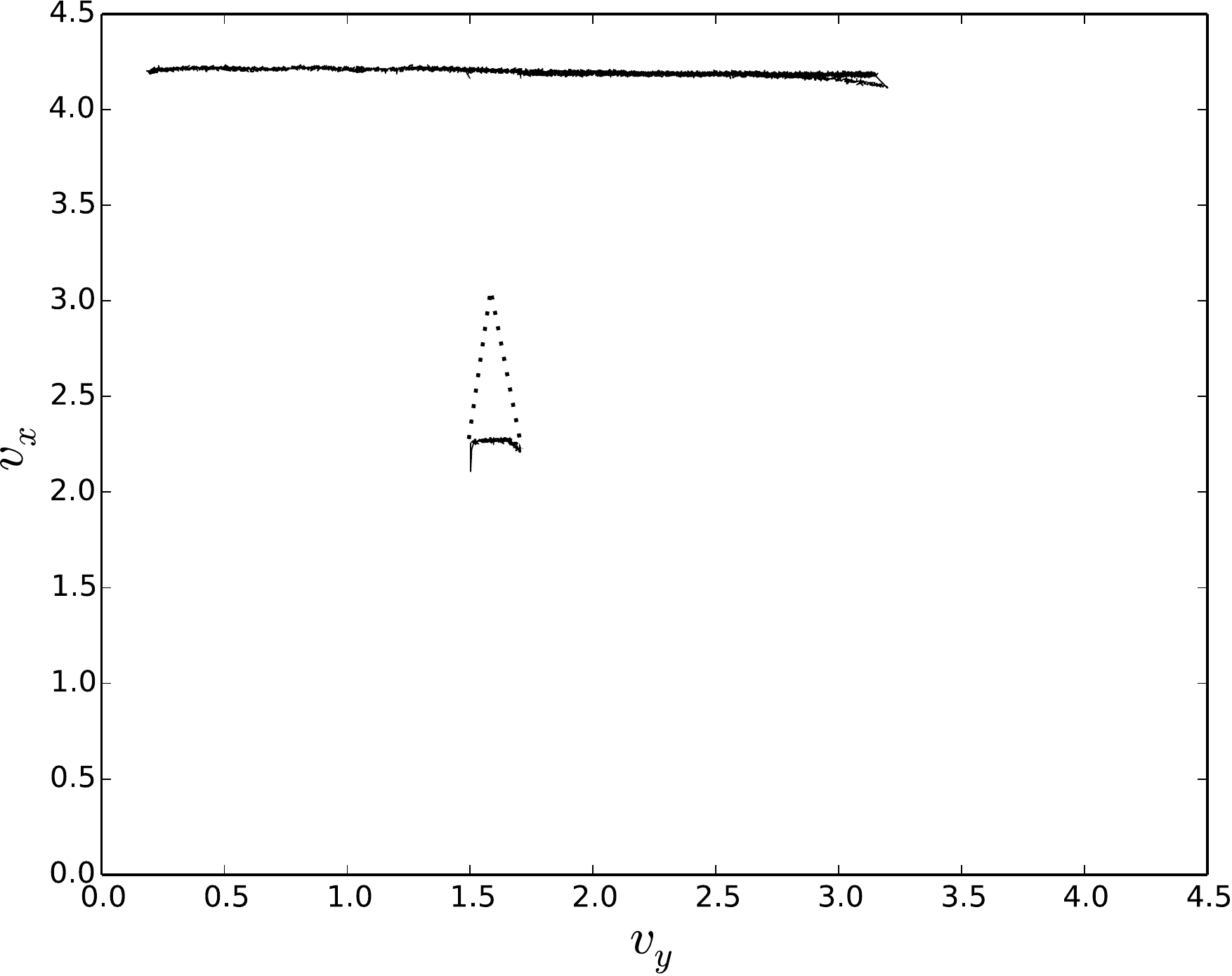}
  \caption{$v_z = 0.8\:V$.}
  \label{fig:exp_mirror_1}
 \end{subfigure}
 \begin{subfigure}[t]{0.32\textwidth}
  \centering
  \includegraphics[scale=0.30]{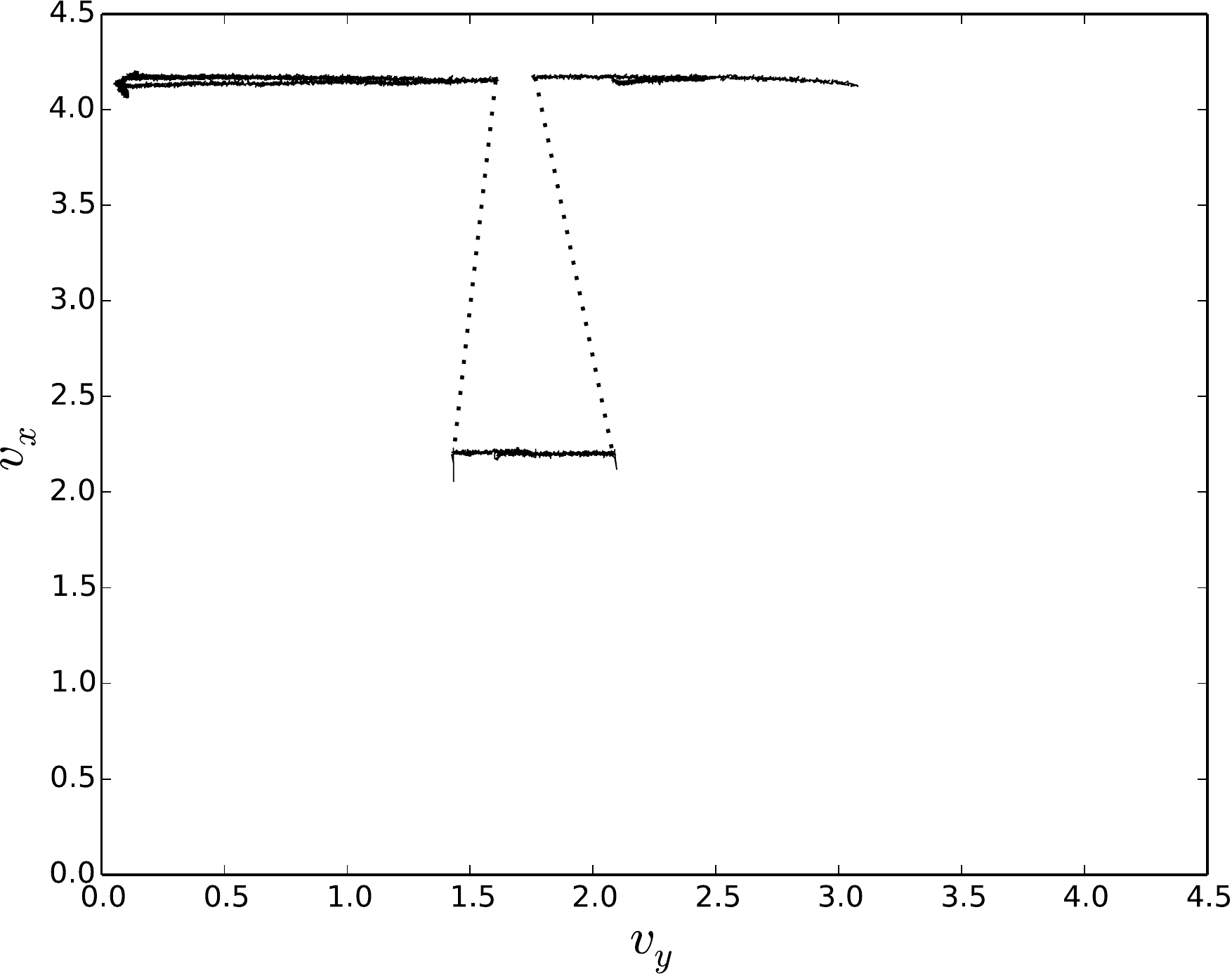}
  \caption{$v_z = 0.9\:V$.}
  \label{fig:exp_mirror_2}
 \end{subfigure}
 \begin{subfigure}[t]{0.32\textwidth}
  \centering
  \includegraphics[scale=0.30]{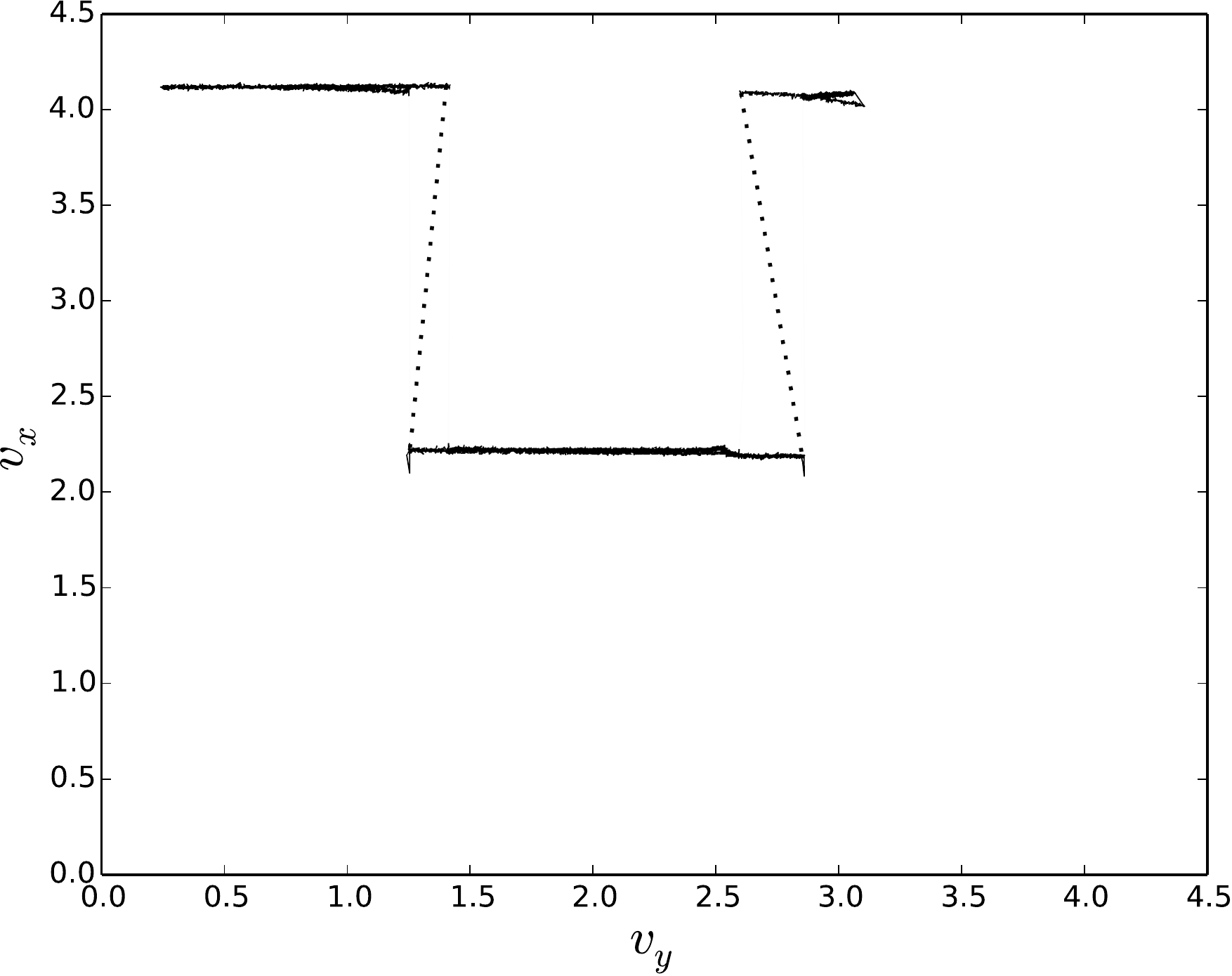}
  \caption{$v_z = 2.7\:V$.}
  \label{fig:exp_mirror_3}
 \end{subfigure}
\caption{Mirrored-hysteresis. Experimental results. Solid lines correspond to the $v_x$ voltage measured while sweeping $v_y$.
 Dotted lines are manually added to sketch the unstable solutions that do not manifest physically
 (cf. Fig.~\ref{fig:mod_mirror_plot}).}
\label{fig:exp_mirror}
\end{figure*}

The resulting circuit is depicted on the upper half of the circuit shown in Fig.~\ref{fig:final_circuit}.
The voltage-controlled non-monotone block is composed of transistors Q1-Q3, the emitter-follower transistor
corresponds to Q4 and the hysteretical block is formed by transistors Q5-Q6. The behaviour is shown in
Fig.~\ref{fig:exp_mirror}. It can be seen that the responses are indeed stiffer than the simulated ones, but
remain qualitatively equivalent to the mirrored hysteresis of Figs.~\ref{FIG:bursting phase portrait}-C. 

The $v_y$--$v_x$ feedback loop is closed by means of another transistor in a common-emitter configuration, Q7,
with a response similar to the one depicted in Fig.~\ref{fig:ce_sat_plot}. The negative slope presented by this block
is corrected by taking the complementary output of $v_x$ at the collector of Q6, instead of $v_x$ itself at the collector of Q5.
The 11\:nF ceramic capacitor and the 100\:k$\Omega$ potentiometer set the slow time-scale. The potentiometer also permits
to adjust the slope of the $v_y$-nullcline in a similar way as in Fig.~\ref{fig:mod_mirror_plot}.
Indeed, by properly adjusting the potentiometer, it is possible to set the circuit either in bursting or in tonic spiking mode.

\begin{figure}
\centering
 \includegraphics[width=0.90\columnwidth]{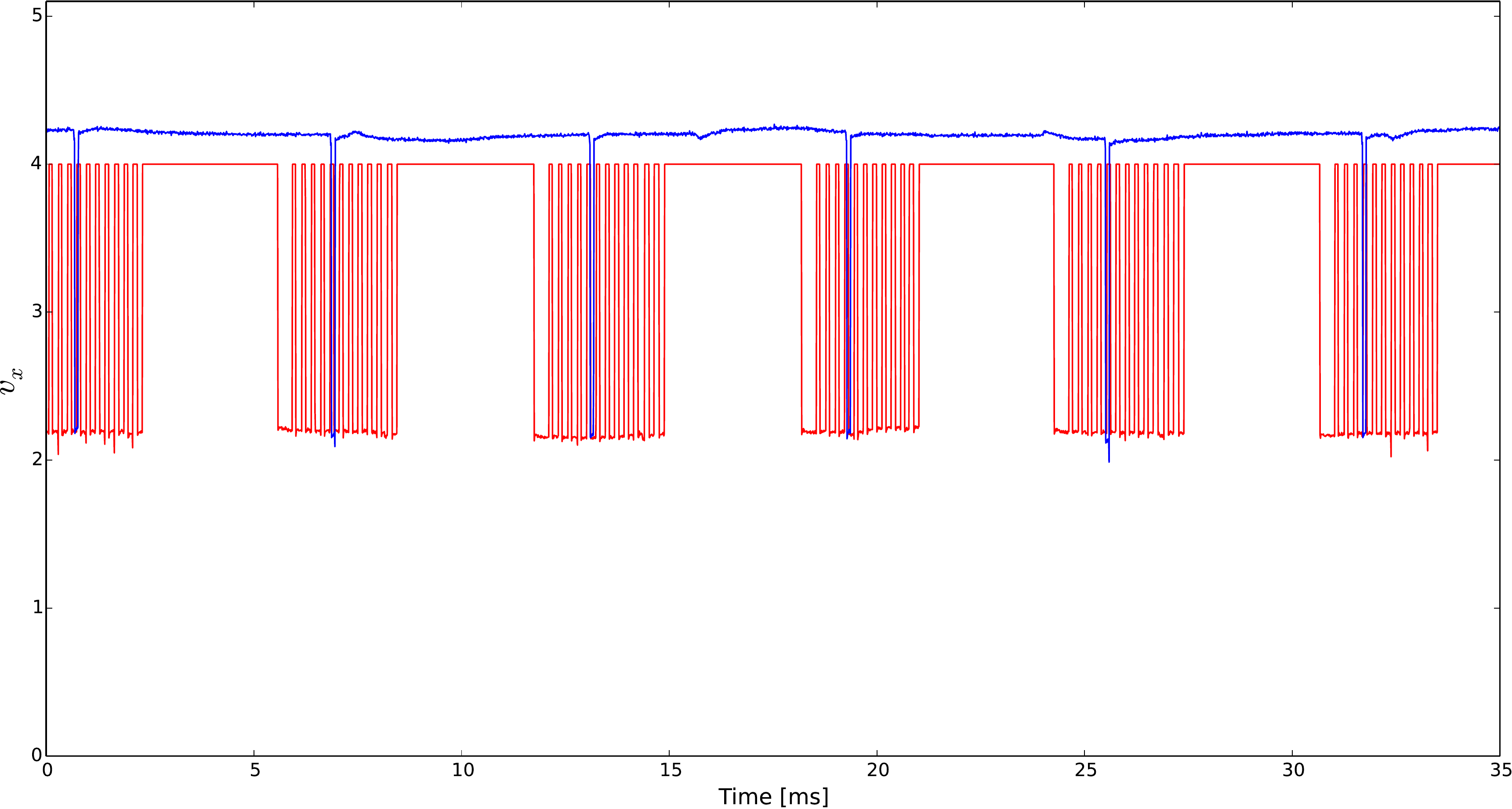}
\caption{Bursting (red) and tonic spiking (blue). Experimental results.}
\label{fig:exp_bs}
\end{figure}

Closing the $v_z$--$v_x$ loop is slightly more complicated. First, an emitter-follower, Q8, is used for the same impedance-matching
purposes as before. A 10\:$\mu$F electrolytic capacitor and a 5.6 \:k$\Omega$ resistor set the ultra slow time-scale. The output
of the ultra slow filter is then connected to a two-stage amplifier, Q9-Q10, that works as a signal conditioner which 
translates the voltage change in $v_x$ that arises when transitioning from the upper stable equilibrium ($\sim$4\:V) to the
stable limit cycle ($\sim$3\:V, average), into the voltage change required in $v_z$ for transitioning between the mirrored
hysteresis of Fig.~\ref{fig:exp_mirror_3} (2.7\:V) and Fig.~\ref{fig:exp_mirror_1} (0.8\:V). The desired bursting and tonic spiking
modes are finally shown in Fig.~\ref{fig:exp_bs}. As explained above, it is possible to smoothly transition from one mode to the other
simply by adjusting the 100\:k$\Omega$ potentiometer.

It is worth mentioning that the circuit was built using common low-cost components. All the transistors have part number 2N2222.
The resistors have precision tolerances of 5\:\%, while the capacitors have precision tolerances of 10\:\%. The low precision of 
the components attests to the intrinsic robust nature of the singularity approach. To further asses the robustness of the design,
the transistors where randomly swapped. The resulting responses were virtually indistinguishable.

\section{Conclusion and perspectives} \label{SEC: conclusions}

{

\subsection{A robust geometric neuromorphic circuit design methodology}

Instead of relying on fine, non-constructive parameter tuning, the design methodology introduced in this paper allows to implement desired behaviors in electronic circuits from the geometrical inspection of the static characteristic of suitable sub-circuits. The proposed methodology is a direct application of the geometric analysis of neuronal behaviors in \cite{Franci2013,Franci2014} and of the realization theory in \cite{Franci2014a}. The main extension was the use of non-smooth analysis to cope with the switching behavior typical of electronic devices.
When applied to the biological transition between tonic spiking and bursting the result of the proposed methodology is a compact six-transistor {\tt ngspice} model that only requires four additional transistors for a robust, low-cost component, real-world implementation.

\subsection{Reliable excitability modulation in neuromorphic circuits}

Because the same geometry of high-dimensional neuron models is being enforced, it is natural to expect that the proposed methodology enforces the same input--output behavior. The transition between tonic spiking and bursting is associated with a switch in excitability type from linear input coding to nonlinear input detection that plays a fundamental role in brain functions. The same qualitative transition is reliably reproduced in the designed neuromorphic devices. Other excitability transitions could robustly and efficiently be implemented following the same geometric methodology, which provides the main advantage with respect to available neuromorphic circuit design methods.

\subsection{Perspective}

The possibility of implementing the transition between distinct excitability types in simple, inexpensive, and robust neuromorphic circuits opens the path, for instance, to the design of neuromorphic sensors inspired by the thalamus, the main sensory hub in the central nervous system, where this transition plays a major role in the efficient coding of sensory stimuli~\cite{Sherman2001}. 

In terms of computational capabilities, spike-based neural networks offer several promising features such as selective attention~\cite{bartolozzi2006} and homeostasis~\cite{bartolozzi2009}.

More generally, because virtually any higher-level brain function relies on the modulation of spiking and the excitability property at the single neuron level \cite{lee2012neuromodulation}, the designed circuit potentially provides a novel building block for any neuromorphic circuit in which neuromodulation is essential.

}

\section*{Acknowledgments}

Alessio Franci acknowledges support by DGAPA-PAPIIT (UNAM) grant IA105816.

\bibliographystyle{IEEEtran}
\bibliography{neuron}

\appendix


Simulations were carried using {\tt ngspice}. The circuit simulator is a well-documented open-source
implementation of {\tt Spice3}, {\tt Cider} and {\tt Xspice}.
The following code describes the circuit depicted in Fig.~\ref{fig:vyz} and was used to produce the plots shown in Fig.~\ref{fig:bs_plot}.

{\small
\begin{verbatim}
* this is complete.cir file
* voltage resources
vcc  5 0 dc 5V
* non monotone
q1  2  3  4 2n2222bis
rC1  5  2 16k
rB1  3  1 100k
rE1  4  0 10k
Ra1  6  1 100k
Ra2  6  2 33k
Rs   6  0 220k
* non monotone modulation
q2  7  8  9 2n2222bis
q3 10 11  9 2n2222bis
rC2  7  5 4.7k
rC3 10  5 4.7k
rB2  8  6 1k
rB3 12 11 1.2k
rE2  9  0 470
* hysteresis
q4  13 14 15 2n2222bis
q5  16 17 15 2n2222bis
rC4 13  5 820
rC5 16  5 240
rB4 14  7 2.4k
rB5 17 13 6k
rE4 15  0 240
* vy-vx feedback loop
ri1 13  1 15k
ri2  1  0 47k 
* Set ri2 to 34.5k for tonic spiking
ciF  1  0 22n
* vz-vx feedback loop
q6  12 18 19 2n2222bis
ro1 16 18 4.7k
ro2 18  0 4.7k
coF 18  0 4.7u 
rC6 12  5 200
rE6 19  0 20
rbi 19  5 150
* model for a 2n2222 transistor
.model 2n2222bis npn (is=14.34f bf=255.9 
+ vaf=74.03 ikf=.2847 ise=14.34f ne=1.307
+ br=6.092 ikr=0 isc=0 nc=2 rb=10 rc=1 
+ cje=22.01p tf=411.1p cjc=7.306p tr=46.91n
+ xtb=1.5 Xti=3 Eg=1.11 Mjc=.3416 Vjc=.75
+ Fc=.5 Mje=.377 Vje=.75 Itf=.6 Vtf=1.7 
+ Xtf=3 )
.control
 tran 1us 40ms
 plot v(16) ylimit 0 5
.endc
.end 
\end{verbatim}
}

\end{document}